\documentclass[11pt]{article}
\usepackage[left=1in,right=1in,top=1in,bottom=1in]{geometry}
\usepackage{times}
\usepackage{expl3}
\usepackage{cite}
\usepackage[table]{xcolor}
\usepackage{multirow}
\usepackage{stackengine} 
\usepackage{hhline}
\usepackage{lipsum}
\usepackage{titlesec}
\usepackage{wrapfig}
\usepackage{epsfig}
\usepackage{graphicx}
\usepackage{amsmath}
\usepackage[title]{appendix}
\usepackage{amssymb}
\usepackage{epstopdf}
\usepackage{boldline}
\usepackage{calligra}
\usepackage{url}
\usepackage{blindtext}

\newcommand{\define}{\stackrel{\mbox{\tiny def}}{=}}

\newtheorem{definition}{Definition}
\newtheorem{theorem}{Theorem}

\newtheorem{lemma}{Lemma}

\usepackage{mathtools}
\usepackage{epstopdf}
\usepackage{balance}
\usepackage{thmtools}
\usepackage{thm-restate}
\usepackage{hyperref}
\usepackage{cleveref}

\usepackage[ruled,vlined]{algorithm2e}
\include{pythonlisting}
\newcommand{\ostar}{\mathbin{\mathpalette\make@circled\star}}

\makeatletter
\newcommand{\removelatexerror}{\let\@latex@error\@gobble}
\makeatother
\makeatletter
\newcommand*{\rom}[1]{\expandafter\@slowromancap\romannumeral #1@}
\makeatother

\ExplSyntaxOn
\newcommand\latinabbrev[1]{
  \peek_meaning:NTF . {
    #1\@}%
  { \peek_catcode:NTF a {
      #1.\@ }%
    {#1.\@}}}
\ExplSyntaxOff

\titleclass{\subsubsubsection}{straight}[\subsubsection]

\begin{document}
\vspace{1cm}
\title{Hanson-Wright Inequality for Random Tensors under Einstein Product}\vspace{1.8cm}
\author{Shih~Yu~Chang 
\thanks{Shih Yu Chang is with the Department of Applied Data Science,
San Jose State University, San Jose, CA, U. S. A. (e-mail: {\tt
shihyu.chang@sjsu.edu}).
           }}

\maketitle

\begin{abstract}
The Hanson-Wright inequality is an upper bound for tails of real quadratic forms in independent subgaussian random variables. In this work, we extend the Hanson-Wright inequality for the maximum eigenvalue of the quadratic sum of random Hermitian tensors under Einstein product. We first prove Weyl inequality for tensors under Einstein product and apply this fact to separate the quadratic form of random Hermitian tensors into diagonal sum and coupling (non-diagonal) sum parts. For the diagonal part, we can apply Bernstein inequality to bound the tail probability of the maximum eigenvalue of the sum of independent random Hermitian tensors directly. For coupling sum part, we have to apply decoupling method first, i.e., decoupling inequality to bound expressions with dependent random Hermitian tensors with independent random Hermitian tensors, before applying Bernstein inequality again to bound the tail probability of the maximum eigenvalue of the coupling sum of independent random Hermitian tensors. Finally, the Hanson-Wright inequality for the maximum eigenvalue of the quadratic sum of random Hermitian tensors under Einstein product can be obtained by the combination of the bound from the diagonal sum part and the bound from the coupling (non-diagonal) sum part. In Appendix of this work, we also include the Hanson-Wright inequality under T-product tensor, which can be derived by the same method of establishing the Hanson-Wright inequality under Einstein product except changing the rule of tensors product operation. 
\end{abstract}


\section{Introduction}\label{sec:Introduction}

The Hanson-Wright inequality provides us an upper bound for tails of real quadratic forms
in independent subgaussian random variables. We define a random variable $X$ is a $\alpha$-subgaussian if for every $\theta > 0$, we have $\mathrm{Pr}( | X | \geq \theta ) \leq 2 \exp ( - \frac{ \theta^2}{2 \alpha^2 })$. The Hanson-Wright inequality states that for any sequence of independent mean zero $\alpha$-subgaussian random variables $X_1, \cdots, X_n$, and any symmetric matrix $\mathbf{A}=(a_{i, j})_{i,j \leq n}$, we have
\begin{eqnarray}\label{eq:HW random variable form}
\mathrm{Pr}\left( \left\vert  \sum\limits_{i,j = 1}^n  a_{i,j}\left(X_i X_j -  \mathbb{E}(X_i X_j    ) \right)                     \right\vert  \geq \theta \right) \leq 2 \exp\left( -\frac{1}{C}\min\left\{ \frac{\theta^2}{\alpha^4 \left\Vert \mathbf{A} \right\Vert_{\mbox{\tiny{HS}}}} , \frac{\theta}{\alpha^2 \left\Vert \mathbf{A} \right\Vert_{\mbox{\tiny{OP}}}} \right\}  \right),
\end{eqnarray}
where $\left\Vert \mathbf{A} \right\Vert_{\mbox{\tiny{HS}}}$ is defined as $\left(  \sum\limits_{i,j = 1}^n  |a_{i,j}|^2   \right)^{1/2}$, and $\left\Vert \mathbf{A} \right\Vert_{\mbox{\tiny{OP}}}$ is defined as $\max\limits_{\left\Vert \mathbf{x}\right\Vert \leq 1} \left\Vert \mathbf{A}\mathbf{x}  \right\Vert_2$. The bound in Eq.~\eqref{eq:HW random variable form} was essentially proved in~\cite{adamczak2015note} in the symmetric case and in~\cite{hanson1971bound} in the zero mean case. The Hanson-Wright inequality has been applied to numerous applications in high-dimensional probability
and statistics, as well as in random matrix theory~\cite{vershynin2018high}. For example, the estimation of bound in Eq.~\eqref{eq:HW random variable form} is applied to the theory of compressed sensing with circulant type matrices~\cite{krahmer2014suprema}. In~\cite{adamczak2015note}, they applied Hanson-Wright inequality to study the concentration properties for sample covariance operators corresponding to Banach space-valued Gaussian random variables.

In recent years, tensors have been applied to different applications in science and engineering~\cite{qi2017tensor}. In data processing fields, tensor theory applications include unsupervised separation of unknown mixtures of data signals~\cite{wu2010robust, mirsamadi2016generalized}, signals filtering~\cite{muti2007survey}, network signal processing~\cite{shen2020topology, shen2017tensor, fu2015joint} and image processing~\cite{ko2020fast, jiang2020framelet}. In wireless communication applications, tensors are applied to model high-dimensional communication channels, e.g., MIMO (multi-input multi-output) code-division~\cite{de2008constrained, zhijin2018blind}, radar communications~\cite{nion2010tensor, sidiropoulos2000parallel}. In numerical multilinear algebra computations, tensors can be applied to solve multilinear system of equations~\cite{wang2019neural}, high-dimensional data fitting/regression~\cite{MR3395816}, tensor complementary problem~\cite{MR3947912}, tensor eigenvalue problem~\cite{MR3479021}, etc. In machine learning, tensors are also used to characterize data with \emph{coupling effects}, for example, tensor decomposition methods have been reported recently to establish the latent-variable models, such as topic models in~\cite{anandkumar2015tensor}, and the method of moments for undertaking the \emph{Latent Dirichlet Allocation} (LDA) in~\cite{sidiropoulos2017tensor}. Nevertheless, all these applications assume that systems modelled by tensors are fixed and such assumption is not true and practical in problems involving tensor formulations. In recent years, there are more works beginning to develop theory about random tensors, see~\cite{MR3616422,MR3783911,MR4140540, chang2021tensor,chang2021TProdI,chang2021TProdII,chang2020convenient,chang2021general}, and references therein. In this work, we extend the Hanson-Wright inequality from random variables to random Hermitian tensors under Einstein product. This is a related work to our another work about Hanson-Wright inequality for symmetric T-product tensor~\cite{HW_T_SYChang_2021}.

We first prove Theorem~\ref{thm:Weyl Inequality for tensor} about Weyl inequality for tensors under Einstein product and apply Theorem~\ref{thm:Weyl Inequality for tensor} to separate the quadratic form of random tensors into diagonal sum and coupling (non-diagonal) sum parts. For the diagonal part, we can apply Theorem~\ref{thm:Bounded Tensor Bernstein_intro} directly to bound the tail probability of the maximum eigenvalue of the sum of independent random Hermitian tensors. For coupling sum part, we have to upper bound this part by Theorem~\ref{thm:decoupling} via decoupling method first, i.e., decoupling inequality to bound expressions with dependent random Hermitian tensors with independent random Hermitian tensors, before applying again Theorem~\ref{thm:Bounded Tensor Bernstein_intro} to bound the tail probability of the maximum eigenvalue of the sum of independent random Hermitian tensors. Therefore, the main result of this work is presented by the following Theorem~\ref{thm:HW inequality}.
\begin{restatable}[Hanson-Wright Inequality for Random Tensors]{thm}{HWThm}\label{thm:HW inequality}
We define a vector of random tensors $\overline{\mathcal{X}} \in \mathbb{R}^{(n \times I_1 \times \cdots \times I_M) \times (I_1 \times \cdots \times I_M)}$ as:
\begin{eqnarray}\label{eq:vec X def 1}
\overline{\mathcal{X}} = \begin{bmatrix}
           \mathcal{X}_{1} \\
           \mathcal{X}_{2} \\
           \vdots \\
           \mathcal{X}_{n}
         \end{bmatrix},
\end{eqnarray}
where random Hermitian tensors $\mathcal{X}_{i} \in \mathbb{R}^{(I_1 \times \cdots \times I_M) \times (I_1 \times \cdots \times I_M)} $ are independent random Hermitian tensors with $\mathbb{E} \mathcal{X}_i = \mathcal{O}$ for $1 \leq i \leq n$. We also require another fixed tensor $\overline{\overline{\mathcal{A}}} \in  \mathbb{R}^{(n \times I_1 \times \cdots \times I_M) \times (n \times I_1 \times \cdots \times I_M)}$, which is defined as:
\begin{eqnarray}\label{eq:matrix A def 1}
\overline{\overline{\mathcal{A}}} = \begin{bmatrix}
           \mathcal{A}_{1,1} &  \mathcal{A}_{1,2} & \cdots & \mathcal{A}_{1, n}  \\
           \mathcal{A}_{2,1} &  \mathcal{A}_{2,2} & \cdots & \mathcal{A}_{2, n}  \\
           \vdots & \vdots &   \vdots & \vdots \\
           \mathcal{A}_{n,1} &  \mathcal{A}_{n,2} & \cdots & \mathcal{A}_{n, n}  \\
         \end{bmatrix},
\end{eqnarray}
where $\mathcal{A}_{i,j} \in \mathbb{R}^{(I_1 \times \cdots \times I_M) \times (I_1 \times \cdots \times I_M)}$ are Hermitian tensors also. We also require following assumptions. Define random Hermitian tensor $\mathcal{Y}_i$ for $i=1,2,\cdots, n$ as 
\begin{eqnarray}\label{eq:diag Yi def main}
\mathcal{Y}_i \define  \mathcal{X}_i \star_M \mathcal{A}_{i,i} \star_M \mathcal{X}_i - \mathbb{E} \left( \mathcal{X}_i   \star_M 
\mathcal{A}_{i,i} \star_M \mathcal{X}_i  \right), 
\end{eqnarray}
we assume that 
\begin{eqnarray}
\mathbb{E} \mathcal{Y}_i = \mathcal{O} \mbox{~~and~~} \lambda_{\max}(\mathcal{Y}_i) \leq T_{\mbox{dg}} 
\mbox{~~almost surely.} 
\end{eqnarray}
Define the total varaince $\sigma_{\mbox{dg}}^2$ as: $\sigma_{\mbox{dg}}^2 \define \left\Vert \sum\limits_{i=1}^n \mathbb{E} \left( \mathcal{Y}^2_i \right) \right\Vert$, where $\left\Vert \cdot \right\Vert$ represents the spectral norm, which equals the largest singular value of a tensor.  The term $\mathbb{I}_1^M$ is defined as the product of each dimension size: $\mathbb{I}_1^M \define \prod\limits_{j=1}^{M} I_j$.

Moreover, we define random Hermitian tensor $\mathcal{Z}_k$ for $k=1,2,\cdots, n^2 -n$ as 
\begin{eqnarray}\label{eq:cp Zk def main}
\mathcal{Z}_k \define   \mathcal{X}^{(1)}_i \star_M  \mathcal{A}_{i,j} \star_M \mathcal{X}^{(2)}_j~~\mbox{for}~1 \leq i \neq j \leq n, 
\end{eqnarray}
where the tensors $\mathcal{X}^{(1)}_i$ are identical distribution copy for the tensors $\mathcal{X}_i$, and the tensors $\mathcal{X}^{(2)}_j$ are identical distribution copy for the tensors $\mathcal{X}_j$, then we assume that 
\begin{eqnarray}
\mathbb{E} \mathcal{Z}_k=\mathcal{O}\mbox{~~almost surely.} 
\end{eqnarray}
Given any realization of the random tensor $\mathcal{X}_i$, denoted as $\tilde{\mathcal{X}}_i$, we assume that 
\begin{eqnarray}\label{eq:max e-value bound cp main}
\lambda_{\max}\left( \tilde{\mathcal{X}}_i \star_M \mathcal{A}_{i, j} \star_M \mathcal{X}_{j} \right) \leq T_{\mbox{cp}} \mbox{~~almost surely},
\end{eqnarray}
where $T_{\mbox{cp}}$ is a positive real number, and all $1 \leq i \neq j \leq n$. The total variance with respect to the tensor $\tilde{\mathcal{X}}_i$ is defined as
\begin{eqnarray}\label{eq:sigma cp main}
\sigma^2_{\mbox{cp}}( \tilde{\mathcal{X}}_i)  \define  \left\Vert \sum\limits_{j=1, \neq i}^{n} \mathbb{E} \left( \left(\tilde{\mathcal{X}}_i \star_M \mathcal{A}_{i, j} \star_M \mathcal{X}_{j}   \right)^2 \right) \right\Vert.
\end{eqnarray}
The function $f(\tilde{\mathcal{X}}_i)$ is the probability density function for the realization tensor $\tilde{\mathcal{X}}_i$.

Then, we have
\begin{eqnarray}
\mathrm{Pr}\left(  \lambda_{\max} \left( \overline{\mathcal{X}}^{\mathrm{T}} 
\overline{\overline{\mathcal{A}}} \overline{\mathcal{X}} - \mathbb{E}\left( \overline{\mathcal{X}}^{\mathrm{T}} 
\overline{\overline{\mathcal{A}}} \overline{\mathcal{X}} \right)  \right) \geq \theta   \right) \leq  
\overbrace{\mathrm{Pr}\left( \lambda_{\max}\left(  \sum\limits_{1 \leq i \neq j \leq n}\mathcal{X}_i \star_M  \mathcal{A}_{i,j} \star_M  \mathcal{X}_j \right) \geq \frac{\theta}{2}  \right) }^{\define \mathrm{P}_{\mbox{cp}} } \nonumber \\
+ \overbrace{\mathrm{Pr} \left( \lambda_{\max}\left(  \sum\limits_{i=1}^{n}  \left(  \mathcal{X}_i \star_M \mathcal{A}_{i,i} \star_M \mathcal{X}_i - \mathbb{E} \left( \mathcal{X}_i  \star_M 
\mathcal{A}_{i,i} \star_M \mathcal{X}_i  \right) \right) \right) \geq \frac{\theta}{2} \right) }^{\define  \mathrm{P}_{\mbox{dg}}} \nonumber \\
 \leq C_4 \mathbb{I}_1^M \sum\limits_{i=1}^{n}  \int_{\tilde{\mathcal{X}}_i}  \exp \left( \frac{ - \theta^2}{  8 n^2C^2_4 \sigma^2_{\mbox{cp}}( \tilde{\mathcal{X}}_i) + 4T_{\mbox{cp}}\theta n C_4 /3 }\right)  f(\tilde{\mathcal{X}}_i) d \tilde{\mathcal{X}}_i  ~~~~~~~~~ \nonumber \\
+  \mathbb{I}_1^M \exp \left( \frac{-\theta^2}{8 \sigma_{\mbox{dg}}^2 + 4 T_{\mbox{dg}}\theta/3}\right),~~~~~~~~~~~~~ ~~~~~~~~~~~~~~~~~~~~~~~~~~~~~~~~~~~~~~~ 
\end{eqnarray}
where $\mathrm{P}_{\mbox{cp}}$ and $\mathrm{P}_{\mbox{dg}}$ are probability bounds related to the coupling sum and the diagonal sum parts, respectively, and the term $C_4$ is a positive constant. 

If $\frac{\theta}{2n C_4} \leq \frac{\sigma^2_{\mbox{cp}}( \tilde{\mathcal{X}}_i)  }{T_{\mbox{cp}}}$ with respect to $i=1,2,\ldots,n$, and $\theta \leq 2 \sigma_{\mbox{dg}}^2 / T_{\mbox{dg}}$, we have 
\begin{eqnarray}
\mathrm{Pr}\left(  \lambda_{\max} \left( \overline{\mathcal{X}}^{\mathrm{T}} 
\overline{\overline{\mathcal{A}}} \overline{\mathcal{X}} - \mathbb{E}\left( \overline{\mathcal{X}}^{\mathrm{T}} 
\overline{\overline{\mathcal{A}}} \overline{\mathcal{X}} \right)  \right) \geq \theta   \right) &\leq&  
C_4 \mathbb{I}_1^M \sum\limits_{i=1}^{n}   \int_{\tilde{\mathcal{X}}_i}  \exp \left( \frac{ - 3 \theta^2}{   32 n^2C^2_4 \sigma^2_{\mbox{cp}}( \tilde{\mathcal{X}}_i) }\right)  f(\tilde{\mathcal{X}}_i) d \tilde{\mathcal{X}}_i  \nonumber \\
&  & + \mathbb{I}_1^M \exp \left( \frac{-3 \theta^2}{ 32 \sigma_{\mbox{dg}}^2}\right). 
\end{eqnarray}
Moreover, if $\frac{\theta}{2n C_4} \geq \frac{\sigma^2_{\mbox{cp}}( \tilde{\mathcal{X}}_i)  }{T_{\mbox{cp}}}$ with respect to $i=1,2,\ldots,n$, and $\theta \geq 2 \sigma_{\mbox{dg}}^2 / T_{\mbox{dg}}$, we have 
\begin{eqnarray}
\mathrm{Pr}\left(  \lambda_{\max} \left( \overline{\mathcal{X}}^{\mathrm{T}} 
\overline{\overline{\mathcal{A}}} \overline{\mathcal{X}} - \mathbb{E}\left( \overline{\mathcal{X}}^{\mathrm{T}} 
\overline{\overline{\mathcal{A}}} \overline{\mathcal{X}} \right)  \right) \geq \theta   \right) &\leq&  n C_4 \mathbb{I}_1^M   \exp \left( \frac{ - 3 \theta}{  16 nC_4 T_{\mbox{cp}}         }\right)   + \mathbb{I}_1^M \exp \left( \frac{-3 \theta}{ 16 T_{\mbox{dg}} } \right). 
\end{eqnarray}
\end{restatable}

The rest of this paper is organized as follows. In Section~\ref{sec:Fundamentals of Tensors and Random Tensors Tail Bounds}, we review tensors under Einstein product and discuss Bernstein bounds for random Hermitian tensors. In Section~\ref{sec:Quadratic Form for Random Tensors and Its Diagonal Sum}, we will formulate our quadratic form for random Hermitan tensors used for Hanson-Wright inequality under Einstein product. Under such quadratic formulation, we will separate this form into diagonal sum and coupling sum parts. The probability bound for the diagonal sum part will also be discussed. The decoupling technique is presented and applied to bound the coupling sum of random Hermitian tensors in Section~\ref{sec:CCoupling Sum of Random Tensors}. Main result of this work: the Hanson-Wright inequality for random Hermitian tensors, is given in Section~\ref{sec:Hanson-Wright Inequality for Random Tensors}. Finally, concluding remarks are given by Section~\ref{sec:Conclusion}.

\section{Fundamentals of Tensors and Random Tensors Tail Bounds}\label{sec:Fundamentals of Tensors and Random Tensors Tail Bounds}

In this section, we will provide a brief introduction of tensors and related theorems in Section~\ref{subsec:Preliminaries of Tensors}. In Section~\ref{subsec:Tensor Bernstein Bounds}, we will present extended Bernstein bounds for a sum of zero-mean random tensors proved in~\cite{chang2020convenient}.

\subsection{Preliminaries of Tensors}\label{subsec:Preliminaries of Tensors}

Throughout this work, scalars are represented by lower-case letters (e.g., $d$, $e$, $f$, $\ldots$), vectors by boldfaced lower-case letters (e.g., $\mathbf{d}$, $\mathbf{e}$, $\mathbf{f}$, $\ldots$), matrices by boldfaced capitalized letters (e.g., $\mathbf{D}$, $\mathbf{E}$, $\mathbf{F}$, $\ldots$), and tensors by calligraphic letters (e.g., $\mathcal{D}$, $\mathcal{E}$, $\mathcal{F}$, $\ldots$), respectively. Tensors are multiarrays of values which are higher-dimensional generalizations from vectors and matrices. Given a positive integer $N$, let $[N] \define \{1, 2, \cdots ,N\}$. An \emph{order-$N$ tensor} (or \emph{$N$-th order tensor}) denoted by $\mathcal{X} \define (a_{i_1, i_2, \cdots, i_N})$, where $1 \leq i_j = 1, 2, \ldots, I_j$ for $j \in [N]$, is a multidimensional array containing $I_1 \times I_2 \times \cdots \times I_{N}$ entries. Let $\mathbb{C}^{I_1 \times \cdots \times I_N}$ and $\mathbb{R}^{I_1 \times \cdots \times I_N}$ be the sets of the order-$N$ $I_1 \times \cdots \times I_N$ tensors over the complex field $\mathbb{C}$ and the real field $\mathbb{R}$, respectively. For example, $\mathcal{X} \in \mathbb{C}^{I_1 \times \cdots \times I_N}$ is an order-$N$ multiarray, where the first, second, ..., and $N$-th dimensions have $I_1$, $I_2$, $\ldots$, and $I_N$ entries, respectively. Thus, each entry of $\mathcal{X}$ can be represented by $a_{i_1, \cdots, i_N}$. For example, when $N = 3$, $\mathcal{X} \in \mathbb{C}^{I_1 \times I_2 \times I_3}$ is a third-order tensor containing entries $a_{i_1, i_2, i_3}$'s.

Without loss of generality, one can partition the dimensions of a tensor into two groups, say $M$ and $N$ dimensions, separately. Thus, for two order-($M$+$N$) tensors: $\mathcal{X} \define (a_{i_1, \cdots, i_M, j_1, \cdots,j_N}) \in \mathbb{C}^{I_1 \times \cdots \times I_M\times
J_1 \times \cdots \times J_N}$ and $\mathcal{Y} \define (b_{i_1, \cdots, i_M, j_1, \cdots,j_N}) \in \mathbb{C}^{I_1 \times \cdots \times I_M\times
J_1 \times \cdots \times J_N}$, according to~\cite{MR3913666}, the \emph{tensor addition} $\mathcal{X} + \mathcal{Y}\in \mathbb{C}^{I_1 \times \cdots \times I_M\times
J_1 \times \cdots \times J_N}$ is given by 
\begin{eqnarray}\label{eq: tensor addition definition}
(\mathcal{X} + \mathcal{Y} )_{i_1, \cdots, i_M, j_1, \cdots,j_N} &\define&
 a_{i_1, \cdots, i_M, j_1 ,\cdots , j_N} \nonumber \\
& &+ b_{i_1, \cdots, i_M, j_1 ,\cdots , j_N}. 
\end{eqnarray}
On the other hand, for tensors $\mathcal{X} \define (a_{i_1, \cdots, i_M, j_1, \cdots,j_N}) \in \mathbb{C}^{I_1 \times \cdots \times I_M\times
J_1 \times \cdots \times J_N}$ and $\mathcal{Y} \define (b_{j_1, \cdots, j_N, k_1, \cdots,k_L}) \in \mathbb{C}^{J_1 \times \cdots \times J_N\times K_1 \times \cdots \times K_L}$, according to~\cite{MR3913666}, the \emph{Einstein product} (or simply referred to as \emph{tensor product} in this work) $\mathcal{X} \star_{N} \mathcal{Y} \in  \mathbb{C}^{I_1 \times \cdots \times I_M\times
K_1 \times \cdots \times K_L}$ is given by 
\begin{eqnarray}\label{eq: Einstein product definition}
\lefteqn{(\mathcal{X} \star_{N} \mathcal{Y} )_{i_1, \cdots, i_M,k_1 ,\cdots , k_L} \define} \nonumber \\ &&\sum\limits_{j_1, \cdots, j_N} a_{i_1, \cdots, i_M, j_1, \cdots,j_N}b_{j_1, \cdots, j_N, k_1, \cdots,k_L}. 
\end{eqnarray}
Note that we will often abbreviate a tensor product $\mathcal{X} \star_{N} \mathcal{Y}$ to ``$\mathcal{X} \hspace{0.05cm}\mathcal{Y}$'' for notational simplicity in the rest of the paper. 
This tensor product will be reduced to the standard matrix multiplication as $L$ $=$ $M$ $=$ $N$ $=$ $1$. Other simplified situations can also be extended as tensor–vector product ($M >1$, $N=1$, and $L=0$) and tensor–matrix product ($M>1$ and $N=L=1$). In analogy to matrix analysis, we define some basic tensors and elementary tensor operations as follows. 

\begin{definition}\label{def: zero tensor}
A tensor whose entries are all zero is called a \emph{zero tensor}, denoted by $\mathcal{O}$. 
\end{definition}

\begin{definition}\label{def: identity tensor}
An \emph{identity tensor} $\mathcal{I} \in  \mathbb{C}^{I_1 \times \cdots \times I_N\times
I_1 \times \cdots \times I_N}$ is defined by 
\begin{eqnarray}\label{eq: identity tensor definition}
(\mathcal{I})_{i_1 , \cdots ,i_N ,
j_1 , \cdots , j_N} \define \prod_{k = 1}^{N} \delta_{i_k, j_k},
\end{eqnarray}
where $\delta_{i_k, j_k} \define 1$ if $i_k  = j_k$; otherwise $\delta_{i_k, j_k} \define 0$.
\end{definition}
In order to define \emph{Hermitian} tensor, the \emph{conjugate transpose operation} (or \emph{Hermitian adjoint}) of a tensor is specified as follows.  
\begin{definition}\label{def: tensor conjugate transpose}
Given a tensor $\mathcal{X} \define (a_{i_1, \cdots, i_M, j_1, \cdots,j_N}) \in \mathbb{C}^{I_1 \times \cdots \times I_M\times J_1 \times \cdots \times J_N}$, its conjugate transpose, denoted by
$\mathcal{X}^{H}$, is defined by
\begin{eqnarray}\label{eq:tensor conjugate transpose definition}
(\mathcal{X}^H)_{ j_1, \cdots,j_N,i_1, \cdots, i_M}  \define  
\overline{a_{i_1, \cdots, i_M,j_1, \cdots,j_N}},
\end{eqnarray}
where the overline notion indicates the complex conjugate of the complex number $a_{i_1, \cdots, i_M,j_1, \cdots,j_N}$. If a tensor $\mathcal{X}$ satisfies $ \mathcal{X}^H = \mathcal{X}$, then $\mathcal{X}$ is a \emph{Hermitian tensor}. 
\end{definition}
\begin{definition}\label{def: unitary tensor}
Given a tensor $\mathcal{U} \define (u_{i_1, \cdots, i_M, j_1, \cdots,j_M}) \in \mathbb{C}^{I_1 \times \cdots \times I_M\times I_1 \times \cdots \times I_M}$, if
\begin{eqnarray}\label{eq:unitary tensor definition}
\mathcal{U}^H \star_M \mathcal{U} = \mathcal{U} \star_M \mathcal{U}^H = \mathcal{I} \in \mathbb{C}^{I_1 \times \cdots \times I_M\times I_1 \times \cdots \times I_M},
\end{eqnarray}
then $\mathcal{U}$ is a \emph{unitary tensor}. 
\end{definition}
In this work, the symbol $\mathcal{U}$ is resrved for a unitary tensor. 

\begin{definition}\label{def: inverse of a tensor}
Given a \emph{square tensor} $\mathcal{Y} \define (a_{i_1, \cdots, i_M, j_1, \cdots,j_M}) \in \mathbb{C}^{I_1 \times \cdots \times I_M\times I_1 \times \cdots \times I_M}$, if there exists $\mathcal{X} \in \mathbb{C}^{I_1 \times \cdots \times I_M\times I_1 \times \cdots \times I_M}$ such that 
\begin{eqnarray}\label{eq:tensor invertible definition}
\mathcal{Y} \star_M \mathcal{X} = \mathcal{X} \star_M \mathcal{Y} = \mathcal{I},
\end{eqnarray}
then $\mathcal{X}$ is the \emph{inverse} of $\mathcal{Y}$. We usually write $\mathcal{X} \define \mathcal{Y}^{-1}$ thereby. 
\end{definition}

We also list other crucial tensor operations here. The \emph{trace} of a square tensor is equivalent to the summation of all diagonal entries such that 
\begin{eqnarray}\label{eq: tensor trace def}
\mathrm{Tr}(\mathcal{X}) \define \sum\limits_{1 \leq i_j \leq I_j,\hspace{0.05cm}j \in [M]} \mathcal{X}_{i_1, \cdots, i_M,i_1, \cdots, i_M}.
\end{eqnarray}
The \emph{inner product} of two tensors $\mathcal{X}$, $\mathcal{Y} \in \mathbb{C}^{I_1 \times \cdots \times I_M\times J_1 \times \cdots \times J_N}$ is given by 
\begin{eqnarray}\label{eq: tensor inner product def}
\langle \mathcal{X}, \mathcal{Y} \rangle \define \mathrm{Tr}\left(\mathcal{X}^H \star_M \mathcal{Y}\right).
\end{eqnarray}
According to Eq.~\eqref{eq: tensor inner product def}, the \emph{Frobenius norm} of a tensor $\mathcal{X}$ is defined by 
\begin{eqnarray}\label{eq:Frobenius norm}
\left\Vert \mathcal{X} \right\Vert \define \sqrt{\langle \mathcal{X}, \mathcal{X} \rangle}.
\end{eqnarray}

We use $\lambda_{\min}$ and $\lambda_{\max}$ to repsent the minimum and the maximum eigenvales of a Hermitain tensor~\cite{ni2019hermitian}. The notation $\succeq$ is used to indicate the semidefinite ordering of tensors. If we have $\mathcal{X} \succeq \mathcal{Y}$, this means that the difference tensor $\mathcal{X} - \mathcal{Y}$ is a positive semidefinite tensor~\cite{ni2019hermitian}.

Following theorem is the min-max theorem for a Hermitian tensor under Einstein product.
\begin{theorem}\label{thm:Courant-Fischer E-product}
Given a Hermitian tensor $\mathcal{C} \in \mathbb{R}^{ I_1 \times \cdots \times I_M \times I_1 \times \cdots \times I_M}$ and $k$ positive integers between $1$ and $~\mathbb{I}_1^M$. Then we have 
\begin{eqnarray}\label{eq:thm:Courant-Fischer E-product}
\lambda_{k} &=_1& \max\limits_{\substack{S \in \mathbb{R}^{ I_1 \times \cdots \times I_M  }\\ \dim(\mathrm{S})  =k  }} \min\limits_{\mathcal{X} \in S } \frac{ \langle \mathcal{X}, \mathcal{C} \star_M \mathcal{X} \rangle }{ \langle \mathcal{X}, \mathcal{X} \rangle }\nonumber \\ 
 &=_2&  \min\limits_{\substack{T \in  \mathbb{R}^{ I_1 \times \cdots \times I_M  }   \\ \dim(T)  = \mathbb{I}_1^M  -k+1   }} \max\limits_{\mathcal{X} \in T} \frac{ \langle \mathcal{X}, \mathcal{C} \star_M \mathcal{X} \rangle }{ \langle \mathcal{X}, \mathcal{X} \rangle },
\end{eqnarray}
where $\lambda_{k}$ is the $k$-th largest eigenvalue of the tensor $\mathcal{C}$. 
\end{theorem}
\textbf{Proof:}
We will just verify the first characterization of $\lambda_{k}$ by $=_1$ in Eq.~\eqref{eq:thm:Courant-Fischer E-product}. The other characterization by $=_2$ in Eq.~\eqref{eq:thm:Courant-Fischer E-product} can be proved similar. For every $\mathcal{X} \in S$ with $\dim(\mathrm{S})  =k$ spanned by $\mathcal{V}_1, \mathcal{V}_2, \cdots, \mathcal{V}_k$ (unitary orthogonal  tensors), we can write $\mathcal{X} = \sum\limits^{k}_{j=1} c_{j}  \mathcal{V}_j$. To show that the value $\lambda_{k}$ is achievable, note that 
\begin{eqnarray}
 \frac{ \langle \mathcal{X}, \mathcal{C} \star \mathcal{X} \rangle }{ \langle \mathcal{X}, \mathcal{X} \rangle }
&=& \frac{  \sum\limits^{k}_{j=1} \lambda_{j}    c_{j}^{\ast} c_{j}          }{    \sum\limits^{k}_{j=1}   c_{j}^{\ast} c_{j}  } \geq  \frac{      \sum\limits^{k}_{j=1} \lambda_{k}  c_{j}^{\ast} c_{j}          }{     \sum\limits^{k}_{j=1}   c_{j}^{\ast} c_{j}  } = \lambda_{k}.
\end{eqnarray}
To verify that this is the maximum, since $T$ has dimension $\mathbb{I}_1^M -k +1$ spanned by $\mathcal{V}_k, \mathcal{V}_{k+1}, \cdots, \mathcal{V}_{\mathbb{I}_1^M }$, we have
\begin{eqnarray}
\min\limits_{\mathcal{X} \in S } \frac{ \langle \mathcal{X}, \mathcal{C} \star \mathcal{X} \rangle }{ \langle \mathcal{X}, \mathcal{X} \rangle } &\leq & \min\limits_{\mathcal{X} \in S \cap T} \frac{ \langle \mathcal{X}, \mathcal{C} \star \mathcal{X} \rangle }{ \langle \mathcal{X}, \mathcal{X} \rangle }.
\end{eqnarray}
Any such $\mathcal{X} \in S \cap T$ can be expressed as $\mathcal{X} =  \sum\limits^{\mathbb{I}_1^M}_{ j=k } c_{j}  \mathcal{V}_{j}$. Then, we have 
\begin{eqnarray}
 \frac{ \langle \mathcal{X}, \mathcal{C} \star \mathcal{X} \rangle }{ \langle \mathcal{X}, \mathcal{X} \rangle }
&=& \frac{  \sum\limits^{\mathbb{I}_1^M}_{j=k} \lambda_{j}    c_{j}^{\ast} c_{j}          }{    \sum\limits^{   \mathbb{I}_1^M  }_{j=k}   c_{j}^{\ast} c_{j}  } \leq  \frac{      \sum\limits^{ \mathbb{I}_1^M }_{j=k} \lambda_{k}  c_{j}^{\ast} c_{j}          }{     \sum\limits^{  \mathbb{I}_1^M  }_{j=k}   c_{j}^{\ast} c_{j}  } = \lambda_{k}.
\end{eqnarray}
Therefore, for all subspaces $S$ of dimensions $k$, we have $\min\limits_{\mathcal{X} \in S} \frac{ \langle \mathcal{X}, \mathcal{C} \star \mathcal{X} \rangle }{ \langle \mathcal{X}, \mathcal{X} \rangle } \leq \lambda_{k}$.
$\hfill \Box$

We then can apply Theorem~\ref{thm:Courant-Fischer E-product} to prove Weyl inequality for Hermitian tensors
under Einstein product.
\begin{theorem}\label{thm:Weyl Inequality for tensor}
Suppose $\mathcal{A}, \mathcal{B} \in \mathbb{C}^{I_1 \times \cdots \times I_M  \times I_1 \times \cdots \times I_M }$ are Hermitian tensors with eigenvalues $\lambda_1 \geq \lambda_2 \geq \cdots \geq \lambda_{\mathbb{I}_1^M}$ and 
$\epsilon_1 \geq \epsilon_2 \geq \cdots \geq \epsilon_{\mathbb{I}_1^M}$, respectively. Let $\mathcal{C} = \mathcal{A} + \mathcal{B}$ with eigenvalues $\mu_1 \geq \mu_2 \geq \cdots \geq \mu_{\mathbb{I}_1^M}$. We then have:
\begin{eqnarray}
\lambda_k + \epsilon_1 \geq_1 \mu_k \geq_2 \lambda_k + \epsilon_{\mathbb{I}_1^M}.
\end{eqnarray}
\end{theorem}
\textbf{Proof:}
Due to Theorem~\ref{thm:Courant-Fischer E-product}, we will prove the inequality $\geq_1$ by $\min\limits_{\substack{T \in  \mathbb{R}^{ I_1 \times \cdots \times I_M  }   \\ \dim(T)  = \mathbb{I}_1^M  -k+1   }} \max\limits_{\mathcal{X} \in T} \frac{ \langle \mathcal{X}, \mathcal{C} \star_M \mathcal{X} \rangle }{ \langle \mathcal{X}, \mathcal{X} \rangle } $ only since the inequality $\geq_2$ can be proved similarly from $ \max\limits_{\substack{S \in \mathbb{R}^{ I_1 \times \cdots \times I_M  }\\ \dim(\mathrm{S})  =k  }} \min\limits_{\mathcal{X} \in S } \frac{ \langle \mathcal{X}, \mathcal{C} \star_M \mathcal{X} \rangle }{ \langle \mathcal{X}, \mathcal{X} \rangle }$.

Because we have
\begin{eqnarray}
\mu_k &=& \min\limits_{\substack{T \in  \mathbb{R}^{ I_1 \times \cdots \times I_M  }   \\ \dim(T)  = \mathbb{I}_1^M  -k+1   }} \max\limits_{\mathcal{X} \in T} \frac{ \langle \mathcal{X}, \left( \mathcal{A} + \mathcal{B}\right)\star_M \mathcal{X} \rangle }{ \langle \mathcal{X}, \mathcal{X} \rangle } \nonumber \\
& = &  \min\limits_{\substack{T \in  \mathbb{R}^{ I_1 \times \cdots \times I_M  }   \\ \dim(T)  = \mathbb{I}_1^M  -k+1   }} \max\limits_{\mathcal{X} \in T} \frac{ \langle \mathcal{X},  \mathcal{A}\star_M \mathcal{X} \rangle +  \langle \mathcal{X},  \mathcal{B}\star_M \mathcal{X} \rangle }{ \langle \mathcal{X}, \mathcal{X} \rangle } \nonumber \\
& \leq & \min\limits_{\substack{T \in  \mathbb{R}^{ I_1 \times \cdots \times I_M  }   \\ \dim(T)  = \mathbb{I}_1^M  -k+1   }} \left( \max\limits_{\mathcal{X} \in T} \frac{ \langle \mathcal{X},  \mathcal{A} \star_M \mathcal{X} \rangle }{ \langle \mathcal{X}, \mathcal{X} \rangle } +
 \max\limits_{\mathcal{X} \in T} \frac{ \langle \mathcal{X},  \mathcal{B} \star_M \mathcal{X} \rangle }{ \langle \mathcal{X}, \mathcal{X} \rangle } \right)  \nonumber \\
& \leq & \min\limits_{\substack{T \in  \mathbb{R}^{ I_1 \times \cdots \times I_M  }   \\ \dim(T)  = \mathbb{I}_1^M  -k+1   }}  \max\limits_{\mathcal{X} \in T} \frac{ \langle \mathcal{X},  \mathcal{A} \star_M \mathcal{X} \rangle }{ \langle \mathcal{X}, \mathcal{X} \rangle } +  \min\limits_{\substack{T \in  \mathbb{R}^{ I_1 \times \cdots \times I_M  }   \\ \dim(T)  = \mathbb{I}_1^M }} 
 \max\limits_{\mathcal{X} \in T} \frac{ \langle \mathcal{X},  \mathcal{B} \star_M \mathcal{X} \rangle }{ \langle \mathcal{X}, \mathcal{X} \rangle } \nonumber \\
&=& \lambda_k + \epsilon_1.
\end{eqnarray}
Then, this theorem is proved. 
$\hfill \Box$

\subsection{Tensor Bernstein Bounds}\label{subsec:Tensor Bernstein Bounds}

Tensor Bernstein inequality is an important inequality to bound the sum of indepdenent, bounded random tensors by restricting the range of the maximum eigenvalue of each random tensors. Following theorem is proved at Theorem 6.2 in~\cite{chang2020convenient}.
\begin{theorem}[Bounded $\lambda_{\max}$ Tensor Bernstein Bounds]\label{thm:Bounded Tensor Bernstein_intro}
Given a finite sequence of independent Hermitian tensors $\{ \mathcal{X}_i  \in \mathbb{C}^{I_1 \times \cdots \times I_M  \times I_1 \times \cdots \times I_M } \}$ that satisfy
\begin{eqnarray}\label{eq1:thm:Bounded Tensor Bernstein_intro}
\mathbb{E} \mathcal{X}_i = 0 \mbox{~~and~~} \lambda_{\max}(\mathcal{X}_i) \leq T 
\mbox{~~almost surely.} 
\end{eqnarray}
Define the total varaince $\sigma^2$ as: $\sigma^2 \define \left\Vert \sum\limits_i^n \mathbb{E} \left( \mathcal{X}^2_i \right) \right\Vert$.
Then, we have following inequalities:
\begin{eqnarray}\label{eq2:thm:Bounded Tensor Bernstein_intro}
\mathrm{Pr} \left( \lambda_{\max}\left( \sum\limits_{i=1}^{n} \mathcal{X}_i \right)\geq \theta \right) \leq \mathbb{I}_1^M \exp \left( \frac{-\theta^2/2}{\sigma^2 + T\theta/3}\right);
\end{eqnarray}
and
\begin{eqnarray}\label{eq3:thm:Bounded Tensor Bernstein_intro}
\mathrm{Pr} \left( \lambda_{\max}\left( \sum\limits_{i=1}^{n} \mathcal{X}_i \right)\geq \theta \right) \leq \mathbb{I}_1^M \exp \left( \frac{-3 \theta^2}{ 8 \sigma^2}\right)~~\mbox{for $\theta \leq \sigma^2/T$};
\end{eqnarray}
and
\begin{eqnarray}\label{eq4:thm:Bounded Tensor Bernstein_intro}
\mathrm{Pr} \left( \lambda_{\max}\left( \sum\limits_{i=1}^{n} \mathcal{X}_i \right)\geq \theta \right) \leq \mathbb{I}_1^M \exp \left( \frac{-3 \theta}{ 8 T } \right)~~\mbox{for $\theta \geq \sigma^2/T$}.
\end{eqnarray}
\end{theorem}

\section{Quadratic Form for Random Hermitian Tensors and Its Diagonal Sum}\label{sec:Quadratic Form for Random Tensors and Its Diagonal Sum}

In Section~\ref{sec:Quadratic Form for Random Tensors}, we will formulate our quadratic form for random tensors used for Hanson-Wright inequality under Einstein product. Under such quadratic formulation, we will separate this form into the diagonal sum and the coupling sum parts. The probability bound for the diagonal sum part will be presented by Section~\ref{sec:Diagonal Sum of Random Tensors}. The coupling sum part will be discussed at next Section~\ref{sec:CCoupling Sum of Random Tensors}.

\subsection{Quadratic Form for Random Tensors}\label{sec:Quadratic Form for Random Tensors}

We define a vector of random tensors $\overline{\mathcal{X}} \in \mathbb{R}^{(n \times I_1 \times \cdots \times I_M) \times (I_1 \times \cdots \times I_M)}$ as:
\begin{eqnarray}\label{eq:vec X def}
\overline{\mathcal{X}} = \begin{bmatrix}
           \mathcal{X}_{1} \\
           \mathcal{X}_{2} \\
           \vdots \\
           \mathcal{X}_{n}
         \end{bmatrix},
\end{eqnarray}
where random Hermitian tensors $\mathcal{X}_{i}$ are independent random tensors with $\mathbb{E} \mathcal{X}_i = \mathcal{O}$ for $1 \leq i \leq n$. We also require another fixed tensor $\overline{\overline{\mathcal{A}}} \in  \mathbb{R}^{(n \times I_1 \times \cdots \times I_M) \times (n \times I_1 \times \cdots \times I_M)}$, which is defined as:
\begin{eqnarray}\label{eq:matrix A def}
\overline{\overline{\mathcal{A}}} = \begin{bmatrix}
           \mathcal{A}_{1,1} &  \mathcal{A}_{1,2} & \cdots & \mathcal{A}_{1, n}  \\
           \mathcal{A}_{2,1} &  \mathcal{A}_{2,2} & \cdots & \mathcal{A}_{2, n}  \\
           \vdots & \vdots &   \vdots & \vdots \\
           \mathcal{A}_{n,1} &  \mathcal{A}_{n,2} & \cdots & \mathcal{A}_{n, n}  \\
         \end{bmatrix},
\end{eqnarray}
where $\mathcal{A}_{i,j} \in \mathbb{R}^{(I_1 \times \cdots \times I_M) \times (I_1 \times \cdots \times I_M)}$ are Hermitian tensors also.

By independence and zero mean of $\mathcal{X}_i$, we can represent $\overline{\mathcal{X}}^{\mathrm{T}} 
\overline{\overline{\mathcal{A}}} \overline{\mathcal{X}} - \mathbb{E}\left( \overline{\mathcal{X}}^{\mathrm{T}} 
\overline{\overline{\mathcal{A}}} \overline{\mathcal{X}} \right)$ as 
\begin{eqnarray}\label{eq:quadratic form diag sum and coupling sum}
\overline{\mathcal{X}}^{\mathrm{T}} 
\overline{\overline{\mathcal{A}}} \overline{\mathcal{X}} - \mathbb{E}\left( \overline{\mathcal{X}}^{\mathrm{T}} 
\overline{\overline{\mathcal{A}}} \overline{\mathcal{X}} \right) &=& \sum\limits_{i=1,j=1}^{n} \mathcal{X}_i \star_M \mathcal{A}_{i,j} \star_M  \mathcal{X}_j  - \sum\limits_{i=1}^{n}  \mathbb{E} \left( \mathcal{X}_i \star_M \mathcal{A}_{i,i} \star_M \mathcal{X}_i \right) \nonumber \\
&=& \sum\limits_{i=1}^{n}\left( \mathcal{X}_i \star_M \mathcal{A}_{i,i} \star_M  \mathcal{X}_i - \mathbb{E} \left( \mathcal{X}_i \star_M  \mathcal{A}_{i,i} \star_M \mathcal{X}_i  \right)\right) \nonumber \\
&  & 
+  \sum\limits_{1 \leq i \neq j \leq n}  \mathcal{A}_{i,j} \star_M \mathcal{X}_i \star_M \mathcal{X}_j
\end{eqnarray}

From Theorem~\ref{thm:Weyl Inequality for tensor} and Eq.~\eqref{eq:quadratic form diag sum and coupling sum}, we have
\begin{eqnarray}
\lambda_{\max} \left( \overline{\mathcal{X}}^{\mathrm{T}} 
\overline{\overline{\mathcal{A}}} \overline{\mathcal{X}} - \mathbb{E}\left( \overline{\mathcal{X}}^{\mathrm{T}} 
\overline{\overline{\mathcal{A}}} \overline{\mathcal{X}} \right)  \right) ~~~~~~~~~~~~~~~~~~~~~~~~~~~~~~~~~~~~~~~~~~~~~~~~~~~~~~~~~~~~~~~~~~~~~~~~~~~ \nonumber \\
 \leq \lambda_{\max}\left(  \sum\limits_{1 \leq i \neq j \leq n} \mathcal{X}_i \star_M  \mathcal{A}_{i,j} \star_M  \mathcal{X}_j \right) + 
\lambda_{\max}\left(  \sum\limits_{i=1}^{n}\left(  \mathcal{X}_i \star_M \mathcal{A}_{i,i} \star_M  \mathcal{X}_i - \mathbb{E} \left( \mathcal{X}_i \star_M  \mathcal{A}_{i,i} \star_M \mathcal{X}_i  \right) \right) \right) 
\end{eqnarray}
Therefore, we have
\begin{eqnarray}\label{eq:prob bpunds by cp and dg}
\mathrm{Pr}\left(  \lambda_{\max} \left( \overline{\mathcal{X}}^{\mathrm{T}} 
\overline{\overline{\mathcal{A}}} \overline{\mathcal{X}} - \mathbb{E}\left( \overline{\mathcal{X}}^{\mathrm{T}} 
\overline{\overline{\mathcal{A}}} \overline{\mathcal{X}} \right)  \right) \geq \theta   \right) \leq  
\overbrace{\mathrm{Pr}\left( \lambda_{\max}\left(  \sum\limits_{1 \leq i \neq j \leq n}\mathcal{X}_i \star_M \mathcal{A}_{i,j} \star_M  \mathcal{X}_j \right) \geq \frac{\theta}{2}  \right) }^{\define \mathrm{P}_{\mbox{cp}} } \nonumber \\
 + \overbrace{\mathrm{Pr} \left( \lambda_{\max}\left(  \sum\limits_{i=1}^n \left(
\mathcal{X}_i \star_M \mathcal{A}_{i,i} \star_M  \mathcal{X}_i - \mathbb{E} \left( \mathcal{X}_i \star_M  \mathcal{A}_{i,i} \star_M \mathcal{X}_i  \right)  \right)  \right) \geq \frac{\theta}{2} \right) }^{\define  \mathrm{P}_{\mbox{dg}}},\end{eqnarray}
where $\mathrm{P}_{\mbox{cp}}$ and $\mathrm{P}_{\mbox{dg}}$ are probability bounds related to the coupling sum and the diagonal sum parts, respectively.

\subsection{Diagonal Sum of Random Tensors}\label{sec:Diagonal Sum of Random Tensors}

The purpose of this subsection is to determine the bound for the probability $\mathrm{P}_{\mbox{dg}}$. If we define the following relation:
\begin{eqnarray}\label{eq:diag Yi def}
\mathcal{Y}_i \define \mathcal{X}_i \star_M \mathcal{A}_{i,i} \star_M  \mathcal{X}_i - \mathbb{E} \left( \mathcal{X}_i \star_M  \mathcal{A}_{i,i} \star_M \mathcal{X}_i  \right) , 
\end{eqnarray}
then, we have $ \mathrm{P}_{\mbox{dg}} $ expressed as:
\begin{eqnarray}
 \mathrm{P}_{\mbox{dg}} = \mathrm{Pr}\left( \lambda_{\max} \left( \sum\limits_{i=1}^n \mathcal{Y}_i \right) \geq \frac{\theta}{2} \right).
\end{eqnarray}

From Theorem~\ref{thm:Bounded Tensor Bernstein_intro}, we will have following lemma about the bound for $ \mathrm{P}_{\mbox{dg}}$.
\begin{lemma}[Bound for $\mathrm{P}_{\mbox{dg}}$]\label{lma:diag sum boud}
Suppose a finite sequence of independent Hermitian tensors \\
$\{ \mathcal{Y}_i  \in \mathbb{C}^{I_1 \times \cdots \times I_M  \times I_1 \times \cdots \times I_M } \}$, which are defined by Eq.~\eqref{eq:diag Yi def},    that satisfy
\begin{eqnarray}\label{eq1:lma:diag sum boud}
\mathbb{E} \mathcal{Y}_i = \mathcal{O} \mbox{~~and~~} \lambda_{\max}(\mathcal{Y}_i) \leq T_{\mbox{dg}} 
\mbox{~~almost surely.} 
\end{eqnarray}
Define the total varaince $\sigma_{\mbox{dg}}^2$ as: $\sigma_{\mbox{dg}}^2 \define \left\Vert \sum\limits_{i=1}^n \mathbb{E} \left( \mathcal{Y}^2_i \right) \right\Vert$.
Then, we have following inequalities:
\begin{eqnarray}\label{eq2:lma:diag sum boud}
\mathrm{P}_{\mbox{dg}}  = \mathrm{Pr} \left( \lambda_{\max}\left( \sum\limits_{i=1}^{n} \mathcal{Y}_i \right)\geq \frac{\theta}{2} \right) \leq \mathbb{I}_1^M \exp \left( \frac{-\theta_{\mbox{dg}} ^2}{ 8 \sigma_{\mbox{dg}}^2 + 4T_{\mbox{dg}} \theta/3}\right);
\end{eqnarray}
and
\begin{eqnarray}\label{eq3:lma:diag sum boud}
\mathrm{P}_{\mbox{dg}}  =  \mathrm{Pr} \left( \lambda_{\max}\left( \sum\limits_{i=1}^{n} \mathcal{Y}_i \right)\geq \frac{\theta}{2} \right) \leq \mathbb{I}_1^M \exp \left( \frac{-3 \theta^2}{ 32 \sigma_{\mbox{dg}}^2}\right)~~\mbox{for $\theta \leq 2\sigma_{\mbox{dg}}^2/T_{\mbox{dg}}$};
\end{eqnarray}
and
\begin{eqnarray}\label{eq4:lma:diag sum boud}
\mathrm{P}_{\mbox{dg}}  =  \mathrm{Pr} \left( \lambda_{\max}\left( \sum\limits_{i=1}^{n} \mathcal{Y}_i \right)\geq \frac{\theta}{2} \right) \leq \mathbb{I}_1^M \exp \left( \frac{-3 \theta}{ 16 T_{\mbox{dg}}} \right)~\mbox{for $\theta \geq 2\sigma_{\mbox{dg}}^2/T_{\mbox{dg}}$}.
\end{eqnarray}
\end{lemma}


\section{Coupling Sum of Random Hermitian Tensors}\label{sec:CCoupling Sum of Random Tensors}

Our next goal is to bound the probability $\mathrm{P}_{\mbox{cp}}$, which is 
\begin{eqnarray}
\mathrm{P}_{\mbox{cp}} &=& \mathrm{Pr}\left( \lambda_{\max}\left(  \sum\limits_{1 \leq i \neq j \leq n} \mathcal{X}_i \star_M  \mathcal{A}_{i,j} \star_M  \mathcal{X}_j \right) \geq \frac{\theta}{2}  \right).
\end{eqnarray}  
However, it is not independent among each summand $\mathcal{X}_i \star_M \mathcal{A}_{i,j} \star_M \mathcal{X}_j $. In order to apply Bernstein bounds in Theorem~\ref{thm:Bounded Tensor Bernstein_intro}, we are interested in doucpling $ \sum\limits_{1 \leq i \neq j \leq n}  \mathcal{X}_i \star_M \mathcal{A}_{i,j} \star_M \mathcal{X}_j $ by
the following expression
\begin{eqnarray}
 \sum\limits_{1 \leq i \neq j \leq n} \mathcal{X}^{(1)}_i \star_M  \mathcal{A}_{i,j} \star_M \mathcal{X}^{(2)}_j,
\end{eqnarray}
where $\{  \mathcal{X}^{(1)}_i\}, \{  \mathcal{X}^{(2)}_i\}$ are two independent copies of random tensors $\{ \mathcal{X}_i \}$. We will prove the following relation:
\begin{eqnarray}\label{eq:decouple rel}
\mathrm{Pr}\left( \lambda_{\max}\left(  \sum\limits_{1 \leq i \neq j \leq n} \mathcal{X}_i \star_M \mathcal{A}_{i,j} \star_M  \mathcal{X}_j \right) \geq \theta  \right) ~~~~~~~~~~~~~~~~~~~~~~~~~~~~~~~~\nonumber \\
\leq  C_4\cdot \mathrm{Pr}\left( \lambda_{\max}\left(   \sum\limits_{1 \leq i \neq j \leq n}  \mathcal{X}^{(1)}_i \star_M  \mathcal{A}_{i,j} \star_M \mathcal{X}^{(2)}_j \right) \geq \frac{\theta}{C_4}  \right),
\end{eqnarray}
where $C_4$ is a positive constant. Our decoupling method discussed in this section is based on the work in~\cite{de1993bounds}, but we extend their approach to the setting of tensors hereof. 

We will present several lemmas before proving the main result of this section. 
\begin{lemma}\label{lma:eq3 eq4}
Let $\mathcal{S}_n$ to be
\begin{eqnarray}\label{eq:eq3}
\mathcal{S}_n & = &  \sum\limits_{1 \leq i \neq j \leq n} \left(\mathcal{X}^{(1)}_i \star_M  \mathcal{A}_{i,j} \star_M  \mathcal{X}^{(1)}_j + \mathcal{X}^{(1)}_i \star_M \mathcal{A}_{i,j} \star_M  \mathcal{X}^{(2)}_j   \right. \nonumber \\
&  & \left.+ \mathcal{X}^{(2)}_i  \star_M \mathcal{A}_{i,j}  \star_M \mathcal{X}^{(1)}_j  +  \mathcal{X}^{(2)}_i \star_M \mathcal{A}_{i,j} \star_M \mathcal{X}^{(2)}_j  \right), 
\end{eqnarray} 
we will have 
\begin{eqnarray}
\mathrm{Pr}\left( \lambda_{\max}\left(  \sum\limits_{1 \leq i \neq j \leq n}  \left(\mathcal{X}^{(1)}_i \star_M  \mathcal{A}_{i,j} \star_M  \mathcal{X}^{(1)}_j   +  \mathcal{X}^{(2)}_i \star_M \mathcal{A}_{i,j} \star_M  \mathcal{X}^{(2)}_j  \right) \right)  \geq  \theta  \right) ~~~~~~  \nonumber \\
\leq \mathrm{Pr}\left(  \lambda_{\max} \left( \mathcal{S}_n \right) \geq \frac{\theta}{3}\right) + 2 \mathrm{Pr} \left( \lambda_{\max} \left(  \sum\limits_{1 \leq i \neq j \leq n}  \left(  \mathcal{X}^{(1)}_i \star_M \mathcal{A}_{i,j} \star_M \mathcal{X}^{(2)}_j    \right) \geq \frac{\theta}{3} \right) \right)
\end{eqnarray} 
\end{lemma}
\textbf{Proof:}
By Theorem~\ref{thm:Weyl Inequality for tensor} and the  triangle inequality of $\lambda_{\max}$. 
$\hfill \Box$

\begin{lemma}\label{lma:lemma1}
Let $\mathcal{X}, \mathcal{Y}$ be two independent and identically distributed random Hermitian tensors with $\mathbb{E}(\mathcal{X}) = \mathbb{E}(\mathcal{Y}) = \mathcal{O}$. Then
\begin{eqnarray}
\mathrm{Pr} \left( \lambda_{\max}\left( \mathcal{X} \right)  \geq \theta \right) \leq 
3 \mathrm{Pr} \left( \lambda_{\max}\left( \mathcal{X} + \mathcal{Y} \right)  \geq \frac{2\theta}{3} \right),
\end{eqnarray}
where $\theta > 0$.
\end{lemma}
\textbf{Proof:}
Let $\mathcal{Z}$ be another independent and identically distributed random Hermitian tensors compared to random Hermitian tensors $\mathcal{X}, \mathcal{Y}$ with $\mathbb{E}(\mathcal{Z})  = \mathcal{O}$. Then, we have
\begin{eqnarray}
\mathrm{Pr} \left( \lambda_{\max}\left( \mathcal{X} \right)  \geq \theta \right)  
&=& \mathrm{Pr} \left( \lambda_{\max}\left( ( \mathcal{X} +  \mathcal{Y}) + 
( \mathcal{X} +  \mathcal{Z}) - ( \mathcal{Y} +  \mathcal{Z})  \right)  \geq 2 \theta \right) \nonumber \\
&\leq& \mathrm{Pr} \left( \lambda_{\max} ( \mathcal{X} +  \mathcal{Y}) \geq \frac{2 \theta}{3}  \right)
+   \mathrm{Pr} \left( \lambda_{\max}\left(  \mathcal{X} +  \mathcal{Z} \right)   \geq \frac{2 \theta}{3}  \right) \nonumber \\
&  & + 
 \mathrm{Pr} \left( \lambda_{\max}\left(  \mathcal{Y} +  \mathcal{Z} \right)   \geq \frac{2 \theta}{3}   \right)
\nonumber \\
&=& 3    \mathrm{Pr} \left( \lambda_{\max}\left(  \mathcal{Y} +  \mathcal{Z} \right)   \geq \frac{2 \theta}{3}   \right)
\end{eqnarray}
$\hfill \Box$

From Lemma~\ref{lma:eq3 eq4}, the proof for the relation provided by Eq.~\eqref{eq:decouple rel} can be reduced as the proof of following two bounds
\begin{eqnarray}\label{eq:5}
\mathrm{Pr}\left( \lambda_{\max}\left(    \sum\limits_{1 \leq i \neq j \leq n}  \mathcal{X}^{(2)}_i \star_M \mathcal{A}_{i,j} \star_M  \mathcal{X}^{(2)}_j    \right) \geq \theta \right) ~~~~~~~~~~~~~~~~~~~~~~~~~~~~~~~~~~~~~~~~~~~~~~~~~~~~~~~~~  \nonumber \\
\leq 3 \mathrm{Pr}\left( \lambda_{\max}\left(    \sum\limits_{1 \leq i \neq j \leq n}  \left( \mathcal{X}^{(1)}_i \star_M  \mathcal{A}_{i,j} \star_M  \mathcal{X}^{(1)}_j  +  \mathcal{X}^{(2)}_i \star_M  \mathcal{A}_{i,j} \star_M  \mathcal{X}^{(2)}_j\right)  \right) \geq \frac{2\theta}{3} \right);
\end{eqnarray}
and, there exists a positive contant $C_2$ to have
\begin{eqnarray}\label{eq:6}
\mathrm{Pr} \left( \lambda_{\max}\left( \mathcal{S}_n \right) \geq \theta \right) \leq C_2
\mathrm{Pr}\left( C_2 \lambda_{\max}\left(    \sum\limits_{1 \leq i \neq j \leq n}   \mathcal{X}^{(1)}_i \star_M \mathcal{A}_{i,j} \star_M  \mathcal{X}^{(2)}_j    \right) \geq \theta \right),
\end{eqnarray}
where the bound provided by Eq.~\eqref{eq:5} is obtained by Lemma~\ref{lma:lemma1}.

We still require two more lemmas before presenting the main result of this section. 
\begin{lemma}\label{lma:prop1}
Let $\mathcal{X} \in \mathfrak{B}$, where $\mathfrak{B}$ is the Banach space with spectral norm,  be any zero mean random Hermitian tensor. Then for all non-random Hermitian tensor $\mathcal{A}$ same dimensions with $\mathcal{X}$ and $ \lambda_{\max}\left( \mathcal{A}\right) > 0$, we have 
\begin{eqnarray}
\mathrm{Pr}\left( \lambda_{\max} \left( \mathcal{A} + \mathcal{X} \right) \geq \lambda_{\max}\left( \mathcal{A}\right) \right) \geq \frac{1}{4} \inf\limits_{f \in F}\frac{  (\mathbb{E}(\left\vert f(\mathcal{X})\right\vert))^2   }{\mathbb{E}(f^2(\mathcal{X}))}
\end{eqnarray}
where $F$ is the family of linear functionals on $\mathfrak{B}$.
\end{lemma}
\textbf{Proof:}
Note that if $x$ is a random variable with $\mathbb{E}x = 0$, then we have $\mathrm{Pr}(x \geq 0) \geq \frac{1}{4} \frac{   (  \mathbb{E}\left\vert x \right\vert )^2 }{\mathbb{E} (x^2)}$. From this fact, we have
\begin{eqnarray}
\mathrm{Pr}\left(f(\mathcal{X}) \geq 0 \right) \geq   \frac{1}{4} \frac{  (\mathbb{E}(\left\vert f(\mathcal{X})\right\vert))^2   }{\mathbb{E}(f^2(\mathcal{X}))},
\end{eqnarray}
since if $f \in F$ is such that $f(\mathcal{A}) = \lambda_{\max}(\mathcal{A})$, then $\{  \lambda_{\max} \left( \mathcal{A} + \mathcal{X} \right) \geq \lambda_{\max}\left( \mathcal{X} \right)   \}$ contains $\{  f \left( \mathcal{A} + \mathcal{X} \right) \geq f \left( \mathcal{A} \right) \} = \{ f \left( \mathcal{X} \right) \geq 0 \}$. 
$\hfill \Box$

\begin{lemma}\label{lma:lemma2}
Let $\mathcal{A}_{i, j}$ for $1 \leq i,j \leq n$, and $\mathcal{B}$ are non-random Hermitian tensors, where $\lambda_{\max}(\mathcal{B}) > 0$. Also let $\{\beta_i \}$ be a sequence of independent and symmetric Bernoulli random variables, that is $\mathrm{Pr}(\beta_i = 1) = \mathrm{Pr}(\beta_i = -1) = \frac{1}{2}$. Then, we have 
\begin{eqnarray}
\mathrm{Pr}\left( \lambda_{\max}\left(  \mathcal{B} + \sum\limits_{i=1}^n  \mathcal{A}_{i, i} \beta_i + \sum\limits_{1 \leq i \neq j, \leq n} \mathcal{A}_{i, j} \beta_i \beta_j \right) \geq \lambda_{\max}(  \mathcal{B} ) \right) \geq C_3
\end{eqnarray}
where $C_3$ is a constant depend on $\left( \sum\limits_{i=1}^n  \mathcal{A}_{i, i} \beta_i + \sum\limits_{1 \leq i \neq j, \leq n} \mathcal{A}_{i, j} \beta_i \beta_j \right) $, but independent of $\mathcal{B}$.
\end{lemma}
\textbf{Proof:}
By setting $\mathcal{X} = \sum\limits_{i=1}^n  \mathcal{A}_{i, i} \beta_i + \sum\limits_{1 \leq i, j, \leq n} \mathcal{A}_{i, j} \beta_i \beta_j $ and $\mathcal{A} = \mathcal{B}$ in Lemma~\ref{lma:prop1}, this lemma is proved. 
$\hfill \Box$

We are ready to present the main Theorem in this section about the bounds on the tail probability by the decoupling inequality.  
\begin{theorem}\label{thm:decoupling}
Given determinstic Hermitian tensors $ \mathcal{A}_{i,j}$ for $1 \leq i, j \leq n$, and random Hermitian tensors $\mathcal{X}_i$ for $1 \leq i \leq n$. Then, there is a positive contant $C_4$ such that for all $n \geq 2$, we have
\begin{eqnarray}
\mathrm{Pr}\left( \lambda_{\max}\left(  \sum\limits_{1 \leq i \neq j \leq n}  \left(\mathcal{X}_i \star_M   \mathcal{A}_{i,j} \star_M  \mathcal{X}_j  \right) \right)  \geq \theta \right) 
\leq ~~~~~~~~~~~~~~~~~~~~~~~~~~~~~~~~~~~~ \nonumber \\ 
C_4 \mathrm{Pr}\left( \lambda_{\max}\left(  \sum\limits_{1 \leq i \neq j \leq n}  \left( \mathcal{X}_i \star_M   \mathcal{A}_{i,j} \star_M  \mathcal{X}^{(2)}_j \right) \right)  \geq \frac{\theta}{C_4} \right),
\end{eqnarray}
where $\theta > 0$. 
\end{theorem}
\textbf{Proof:}
From Lemma~\ref{lma:eq3 eq4}, this theorem can be proved by proving following two bounds:
\begin{eqnarray}\label{eq:eq5}
\mathrm{Pr}\left( \lambda_{\max}\left(  \sum\limits_{1 \leq i \neq j \leq n}  \left( \mathcal{X}^{(2)}_i \star_M  \mathcal{A}_{i,j} \star_M  \mathcal{X}^{(2)}_j  \right) \right)  \geq \theta \right) ~~~~~~~~~~~~~~~~~~~~~~~~~~~~~~~~~~~~~~~~~~~~~~~~~~~~ \nonumber \\
\leq 3 \mathrm{Pr}\left( \lambda_{\max}\left(  \sum\limits_{1 \leq i \neq j \leq n}  \left( \mathcal{X}^{(1)}_i \star_M   \mathcal{A}_{i,j} \star_M \mathcal{X}^{(1)}_j   +  \mathcal{X}^{(2)}_i \star_M \mathcal{A}_{i,j} \star_M \mathcal{X}^{(2)}_j \right) \right)  \geq \frac{2 \theta}{3} \right),
\end{eqnarray}
and
\begin{eqnarray}\label{eq:eq6}
\mathrm{Pr}\left(  \lambda_{\max}\left( \mathcal{S}_n \right) \geq \theta \right) \leq C_2
\mathrm{Pr}\left(C_2  \lambda_{\max}\left(  \sum\limits_{1 \leq i \neq j \leq n}  \left(  \mathcal{X}^{(1)}_i \star_M  \mathcal{A}_{i,j} \star_M \mathcal{X}^{(2)}_j  \right) \right)  \geq \theta \right) ,
\end{eqnarray}
where $C_2$ is a constant.

To prove Eq.~\eqref{eq:eq6}, we first transform the problem of proving Eq.~\eqref{eq:eq6} into a problem
conditionally with a non-homogeneous binomial in Bernoulli random variables. Let $\{\rho_i \}$ be a sequence of independent and symmetric Bernoulli random variables independent of random Hermitian tensors $\{\mathcal{X}^{(1)}_i\}, \{\mathcal{X}^{(2)}_i\}$. Let $(\mathcal{Z}^{(1)}, \mathcal{Z}^{(2)}) = (\mathcal{X}^{(1)}, \mathcal{X}^{(2)}) $ if $\rho_i = 1$, and  $(\mathcal{Z}^{(1)}, \mathcal{Z}^{(2)}) = (\mathcal{X}^{(2)}, \mathcal{X}^{(1)}) $ if $\rho_i = -1$. Then, we have
\begin{eqnarray}\label{eq:eq8}
4  \mathcal{Z}^{(1)}_i \star_M \mathcal{A}_{i,j} \star_M \mathcal{Z}^{(2)}_j = 
~~~~~~~~~~~~~~~~~~~~~~~~~~~~~~~~~~~~~~~~~~~~~~~~~~~~~~~~~~~~~~~~~~~~~~~~~~~~~~~~~~~~~~~~~~~~~~~~~~~~~~~\nonumber \\
(1 + \rho_i) (1 - \rho_j)   \mathcal{X}^{(1)}_i \star_M \mathcal{A}_{i,j} \star_M  \mathcal{X}^{(1)}_j+
(1 + \rho_i) (1 + \rho_j)  \mathcal{X}^{(1)}_i \star_M \mathcal{A}_{i,j} \star_M  \mathcal{X}^{(2)}_j \nonumber \\
+ (1 - \rho_i) (1 - \rho_j) \mathcal{X}^{(2)}_i \star_M \mathcal{A}_{i,j} \star_M \mathcal{X}^{(1)}_j  +
(1 - \rho_i) (1 + \rho_j)   \mathcal{X}^{(2)}_i \star_M \mathcal{A}_{i,j} \star_M  \mathcal{X}^{(2)}_j.
\end{eqnarray}
If we define $\mathfrak{P}$ as a realization of $\rho_i$ for $1 \leq i \leq n$, we have
\begin{eqnarray}\label{eq:eq9}
4 \mathbb{E}( \mathcal{Z}^{(1)}_i \star_M \mathcal{A}_{i,j} \star_M  \mathcal{Z}^{(2)}_j   | \mathfrak{P} )
&=& \{ \mathcal{X}^{(1)}_i \star_M  \mathcal{A}_{i,j} \star_M  \mathcal{X}^{(1)}_j + 
        \mathcal{X}^{(1)}_i \star_M \mathcal{A}_{i,j} \star_M  \mathcal{X}^{(2)}_j  \nonumber \\
&  &+ 
        \mathcal{X}^{(2)}_i \star_M \mathcal{A}_{i,j} \star_M  \mathcal{X}^{(1)}_j +
         \mathcal{X}^{(2)}_i \star_M \mathcal{A}_{i,j} \star_M \mathcal{X}^{(2)}_j\}.
\end{eqnarray}

By setting $\mathcal{B} = \mathcal{S}_n$ in Lemma~\ref{lma:lemma2}, and Eqs.~\eqref{eq:eq3}~\eqref{eq:eq8},~\eqref{eq:eq9}, we have
\begin{eqnarray}\label{eq:eq9-1}
\mathrm{Pr}\left( 4 \lambda_{\max}\left(    \sum\limits_{1 \leq i \neq j \leq n}  \left(\mathcal{Z}^{(1)}_i \star_M   \mathcal{A}_{i,j} \star_M  \mathcal{Z}^{(2)}_j      \right)  \right)  \geq \lambda_{\max}(\mathcal{S}_n)  \Bigg\vert \mathfrak{P} \right) \geq C_4.
\end{eqnarray}
By integrating over the set $\{  \lambda_{\max}(\mathcal{S}_n) \geq \theta \}$, we obtain
\begin{eqnarray}\label{eq:thm:decoupling}
\frac{1}{C_4} \mathrm{Pr}\left( \lambda_{\max}( \mathcal{S}_n )  \geq \theta \right) & \leq &
\mathrm{Pr}\left(4 
 \lambda_{\max} \left(     \sum\limits_{1 \leq i \neq j \leq n}  \left( \mathcal{Z}^{(1)}_i \star_M  \mathcal{A}_{i,j} \star_M  \mathcal{Z}^{(2)}_j          \right)\right)  \geq \theta \right) \nonumber \\
&=& \mathrm{Pr}\left( 4 
 \lambda_{\max} \left(     \sum\limits_{1 \leq i \neq j \leq n}  \left( \mathcal{X}^{(1)}_i \star_M  \mathcal{A}_{i,j} \star_M  \mathcal{X}^{(2)}_j          \right)\right)  \geq \theta \right),
\end{eqnarray}
because the sequence $\{  \mathcal{X}^{(1)}_i, \mathcal{X}^{(2)}_i \}$ for $1 \leq i \leq n$ has the same distribution as  $\{  \mathcal{Z}^{(1)}_i, \mathcal{Z}^{(2)}_i \}$ for $1 \leq i \leq n$. The proof is completed 
by using inequality in Eq.~\eqref{eq:thm:decoupling} with inequalities in Lemma~\ref{lma:eq3 eq4} and Eq.~\eqref{eq:eq5}.
$\hfill \Box$

We are ready to determine the bound for the probability $\mathrm{P}_{\mbox{cp}}$. If we define the following relation:
\begin{eqnarray}\label{eq:cp Zk def}
\mathcal{Z}_k \define \mathcal{X}^{(1)}_i \star_M  \mathcal{A}_{i,j} \star_M  \mathcal{X}^{(2)}_j,
\end{eqnarray}
then, we have $ \mathrm{P}_{\mbox{cp}} $ expressed as:
\begin{eqnarray}\label{eq:Bound for Pcp-1}
\mathrm{P}_{\mbox{cp}} &=& \mathrm{Pr}\left( \lambda_{\max}\left(  \sum\limits_{1 \leq i \neq j \leq n}\mathcal{X}_i \star_M   \mathcal{A}_{i,j} \star_M \mathcal{X}_j \right) \geq \frac{\theta}{2}  \right) \nonumber \\
&\leq_1& C_4 \mathrm{Pr}\left( \lambda_{\max}\left(  \sum\limits_{k=1}^{n^2 - n} \mathcal{Z}_k  \right)  \geq \frac{\theta}{2C_4} \right),
\end{eqnarray}
where $\leq_1$ is due to Theorem~\ref{thm:decoupling}.

From Theorem~\ref{thm:decoupling}, we will have following lemma about the bound for $ \mathrm{P}_{\mbox{cp}}$.
\begin{lemma}[Bound for $\mathrm{P}_{\mbox{cp}}$]\label{lma:cp sum boud}
Given any realization of the random tensor $\mathcal{X}_i$, denoted as $\tilde{\mathcal{X}}_i$, we assume that 
\begin{eqnarray}\label{eq:max e-value bound cp}
\lambda_{\max}\left( \tilde{\mathcal{X}}_i \star_M \mathcal{A}_{i, j} \star_M \mathcal{X}_{j} \right) \leq T_{\mbox{cp}} \mbox{~~almost surely},
\end{eqnarray}
where $T_{\mbox{cp}}$ is a positive real number, and all $i,j \in \{1,2,\ldots, n\}$ with $i \neq j$. The total variance with respect to the tensor $\tilde{\mathcal{X}}_i$ is defined as
\begin{eqnarray}\label{eq:sigma cp}
\sigma^2_{\mbox{cp}}( \tilde{\mathcal{X}}_i)  \define  \left\Vert \sum\limits_{j=1, \neq i}^{n} \mathbb{E} \left( \left(\tilde{\mathcal{X}}_i \star_M \mathcal{A}_{i, j} \star_M \mathcal{X}_{j}   \right)^2 \right) \right\Vert.
\end{eqnarray}
The function $f(\tilde{\mathcal{X}}_i)$ is the probability density function for the realization tensor $\tilde{\mathcal{X}}_i$. We also assume that $\mathcal{X}_i  \sum\limits_{j=1, \neq i }^{n}\mathcal{A}_{i, j} \mathcal{X}_j $ are Hermitian tensors for $i=1,2,\ldots,n$.

Then, we have following inequalities~\footnote{Note that all superscripts $\mathcal{X}^{(1)}_i$ are removed since they are random copies of $\mathcal{X}_i$.}:
\begin{eqnarray}\label{eq2:lma:cp sum boud}
\mathrm{P}_{\mbox{cp}}  \leq C_4 \mathrm{Pr} \left( \lambda_{\max}\left( \sum\limits_{k=1}^{n^{2} - n} \mathcal{Z}_k \right)\geq \frac{\theta}{2 C_4 } \right) \leq ~~~~~~~~~~~~~~~~~~~~~~~~~~~~~~~~~~~~~~~ \nonumber \\
C_4 \mathbb{I}_1^M \sum\limits_{i=1}^{n}  \int_{\tilde{\mathcal{X}}_i}  \exp \left( \frac{ - \theta^2}{  8 n^2C^2_4 \sigma^2_{\mbox{cp}}( \tilde{\mathcal{X}}_i) + 4T_{\mbox{cp}}\theta n C_4 /3 }\right)  f(\tilde{\mathcal{X}}_i) d \tilde{\mathcal{X}}_i ;
\end{eqnarray}
and
\begin{eqnarray}\label{eq3:lma:cp sum boud}
\mathrm{P}_{\mbox{cp}}  \leq  C_4  \mathrm{Pr} \left( \lambda_{\max}\left( \sum\limits_{k=1}^{n^2 - n} \mathcal{Z}_k \right)\geq \frac{\theta}{2 C_4 } \right) & \leq & 
 C_4 \mathbb{I}_1^M \sum\limits_{i=1}^{n}   \int_{\tilde{\mathcal{X}}_i}  \exp \left( \frac{ - 3 \theta^2}{   32 n^2C^2_4 \sigma^2_{\mbox{cp}}( \tilde{\mathcal{X}}_i) }\right)  f(\tilde{\mathcal{X}}_i) d \tilde{\mathcal{X}}_i \nonumber \\ 
&  & \mbox{for $\frac{\theta}{2n C_4} \leq \frac{\sigma^2_{\mbox{cp}}( \tilde{\mathcal{X}}_i)  }{T_{\mbox{cp}}}$ with respect to $i=1,\ldots,n$;}
\end{eqnarray}
and
\begin{eqnarray}\label{eq4:lma:cp sum boud}
\mathrm{P}_{\mbox{cp}}  \leq  \mathrm{Pr} \left( \lambda_{\max}\left( \sum\limits_{k=1}^{n^2 - n} \mathcal{Z}_k \right)\geq \frac{\theta}{2 C_4} \right) &\leq&  n  C_4 \mathbb{I}_1^M 
  \exp \left( \frac{ - 3 \theta}{  16 nC_4 T_{\mbox{cp}}         }\right) 
\nonumber \\
& &\mbox{for  $\frac{\theta}{2n C_4} \geq \frac{\sigma^2_{\mbox{cp}}( \tilde{\mathcal{X}}_i)  }{T_{\mbox{cp}}}$ with respect to $i=1,2,\ldots,n$.}
\end{eqnarray}
\end{lemma}
\textbf{Proof:}
Since all Einstein product are same in this proof, we will remove $\star_M$ for space saving in this proof. From Eq.~\eqref{eq:Bound for Pcp-1}, we have 
\begin{eqnarray}\label{eq5:lma:cp sum boud}
C_4 \mathrm{Pr}\left( \lambda_{\max}\left(  \sum\limits_{k=1}^{n^2 - n} \mathcal{Z}_k  \right)  \geq \frac{\theta}{2C_4} \right) = 
C_4 \mathrm{Pr}\left( \lambda_{\max}\left(  \sum\limits_{i=1}^{n}\left( \mathcal{X}^{(1)}_i  \sum\limits_{j=1, \neq i }^{n}\mathcal{A}_{i, j} \mathcal{X}^{(2)}_j \right)  \right)  \geq \frac{\theta}{2C_4} \right) \nonumber \\
\leq_1 C_4  \sum\limits_{i=1}^{n} \mathrm{Pr}\left( \lambda_{\max} \left(  \mathcal{X}^{(1)}_i \left( \sum\limits_{j=1, \neq i }^{n}\mathcal{A}_{i, j} \mathcal{X}^{(2)}_j  \right) \right)  \geq \frac{\theta}{2n C_4}\right),
\end{eqnarray}
where we apply Theorem~\ref{thm:Weyl Inequality for tensor} in $\leq_1$. By conditional probability with respect to $\mathcal{X}_i$, each term in Eq.~\eqref{eq5:lma:cp sum boud} can be expressed as 
\begin{eqnarray}\label{eq6:lma:cp sum boud}
\mathrm{Pr}\left( \lambda_{\max} \left(  \mathcal{X}^{(1)}_i \left( \sum\limits_{j=1, \neq i }^{n}\mathcal{A}_{i, j} \mathcal{X}^{(2)}_j  \right) \right)  \geq \frac{\theta}{2n C_4}\right) = ~~~~~~~~~~~~~~~~~~~~~~~~~~~~~~~~~~~~~~~~~~~~~~~ \nonumber \\
\int_{\tilde{\mathcal{X}}^{(1)}_i} \mathrm{Pr}\left( \lambda_{\max} \left(  \tilde{\mathcal{X}}^{(1)}_i \left( \sum\limits_{j=1, \neq i }^{n}\mathcal{A}_{i, j} \mathcal{X}^{(2)}_j  \right) \right)  \geq \frac{\theta}{2n C_4}\right) f(\tilde{\mathcal{X}}^{(1)}_i) d \tilde{\mathcal{X}}^{(1)}_i.
\end{eqnarray}

From Theorem~\ref{thm:Bounded Tensor Bernstein_intro} and Eq.~\eqref{eq6:lma:cp sum boud} with conditions given by Eqs.~\eqref{eq:max e-value bound cp} and~\eqref{eq:sigma cp}, we have
\begin{eqnarray}\label{eq7:lma:cp sum boud}
\mathrm{Pr}\left( \lambda_{\max} \left(  \mathcal{X}^{(1)}_i \left( \sum\limits_{j=1, \neq i }^{n}\mathcal{A}_{i, j} \mathcal{X}^{(2)}_j  \right) \right)  \geq \frac{\theta}{2n C_4}\right) \leq  ~~~~~~~~~~~~~~~~~~~~~~~~~~~~~~~~~~~~~~~~~~~~~~~ \nonumber \\
\mathbb{I}_1^M \int_{\tilde{\mathcal{X}}^{(1)}_i}  \exp \left( \frac{ - \theta^2}{  8 n^2C^2_4 \sigma^2_{\mbox{cp}}( \tilde{\mathcal{X}}^{(1)}_i) + 4T_{\mbox{cp}}\theta n C_4 /3 }\right)  f(\tilde{\mathcal{X}}^{(1)}_i) d \tilde{\mathcal{X}}^{(1)}_i.
\end{eqnarray}
If $\frac{\theta}{2n C_4} \leq \frac{\sigma^2_{\mbox{cp}}( \tilde{\mathcal{X}}_i)  }{T_{\mbox{cp}}}$ for all $i = 1,2,\ldots, n$, we can have the following bound 
\begin{eqnarray}\label{eq8:lma:cp sum boud}
\mathrm{Pr}\left( \lambda_{\max} \left(  \mathcal{X}^{(1)}_i \left( \sum\limits_{j=1, \neq i }^{n}\mathcal{A}_{i, j} \mathcal{X}^{(2)}_j  \right) \right)  \geq \frac{\theta}{2n C_4}\right) \leq  ~~~~~~~~~~~~~~~~~~~~~~~~~~~~~~~~~~~~~~~~~~~~~~~ \nonumber \\
\mathbb{I}_1^M \int_{\tilde{\mathcal{X}}^{(1)}_i}  \exp \left( \frac{ - 3 \theta^2}{   32 n^2C^2_4 \sigma^2_{\mbox{cp}}( \tilde{\mathcal{X}}^{(1)}_i) }\right)  f(\tilde{\mathcal{X}}^{(1)}_i) d \tilde{\mathcal{X}}^{(1)}_i.
\end{eqnarray}
On the other hand, if $\frac{\theta}{2n C_4} \geq \frac{\sigma^2_{\mbox{cp}}( \tilde{\mathcal{X}}_i)  }{T_{\mbox{cp}}}$ for all $i = 1,2,\ldots, n$, we can have the following bound 
\begin{eqnarray}\label{eq9:lma:cp sum boud}
\mathrm{Pr}\left( \lambda_{\max} \left(  \mathcal{X}^{(1)}_i \left( \sum\limits_{j=1, \neq i }^{n}\mathcal{A}_{i, j} \mathcal{X}^{(2)}_j  \right) \right)  \geq \frac{\theta}{2n C_4}\right) \leq  ~~~~~~~~~~~~~~~~~~~~~~~~~~~~~~~~~~~~~~~~~~~~~~~ \nonumber \\
\mathbb{I}_1^M \int_{\tilde{\mathcal{X}}^{(1)}_i}  \exp \left( \frac{ - 3 \theta}{  16 nC_4 T_{\mbox{cp}}         }\right)  f(\tilde{\mathcal{X}}^{(1)}_i) d \tilde{\mathcal{X}}^{(1)}_i =
\mathbb{I}_1^M  \exp \left( \frac{ - 3 \theta}{  16 nC_4 T_{\mbox{cp}}         }\right). 
\end{eqnarray}
This Lemma is proved by combining Eq.~\eqref{eq5:lma:cp sum boud} with Eqs.~\eqref{eq7:lma:cp sum boud}~\eqref{eq8:lma:cp sum boud} and~\eqref{eq9:lma:cp sum boud}, respectively.
$\hfill \Box$

\section{Hanson-Wright Inequality for Random Tensors}\label{sec:Hanson-Wright Inequality for Random Tensors}

In this section, we will present the proof for the main result of this paper, the Hanson-Wright inequality for random Hermitian tensors. 

%
\textbf{Proof:}
By combining Lemma~\ref{lma:diag sum boud} and Lemma~\ref{lma:cp sum boud} with Eq.~\eqref{eq:prob bpunds by cp and dg}, this theorem is proved.
$\hfill \Box$

\section{Conclusion}\label{sec:Conclusion} 

In this work, we generalize the Hanson-Wright inequality from the quadratic forms
in independent subgaussian random variables to the random Hermitian tensors. First, we apply Weyl inequality for tensors under the Einstein product and apply this fact to separate the quadratic form of random Hermitian tensors into the diagonal sum and the coupling (non-diagonal) sum parts. Second, we apply decoupling inequality to bound expressions with dependent random Hermitian tensors with independent random Hermitian tensors. Finally, the Hanson-Wright inequality can be obtained by utilizing Bernstein inequality to the diagonal sum part and the coupling sum part, respectively.  \\\\

\textbf{\Large Appendix:  Hanson-Wright inequality for T-product tensors}\\

\begin{appendices}

The T-product operation between two three order tensors was introduced by Kilmer and her collaborators in~\cite{kilmer2013third}. In this Appendix, we will apply the same technique used in the previous sections to first establish Courant-Fischer theorem for a T-product tensor in Appendix~\ref{sec:Courant-Fischer Theorem for T-product tensor}, and use this fact to build Weyl inequality for symmetric T-product tensors in Appendix~\ref{sec:Weyl Inequality for Symmetric T-product Tensors}. Finally, we have the Hanson-Wright inequality for random symmetric T-product tensors presented by Appendix~\ref{sec:Hanson-Wright Inequality for Random Symmetric T-product Tensors}.

\section{Courant-Fischer Theorem for T-product tensor}\label{sec:Courant-Fischer Theorem for T-product tensor}

If a T-product tensor $\mathcal{C} \in \mathbb{R}^{m \times m \times p}$ can be diagonalized as
\begin{eqnarray}\label{eq:block diagonalized format}
\mbox{bcirc}(\mathcal{C}) = \left( \mathbf{F}^{\mathrm{H}}_m \otimes \mathbf{I}_m \right) \mbox{Diag}\left( \mathbf{C}_i: i \in \{1, \cdots, m \} \right)  \left( \mathbf{F}_m \otimes \mathbf{I}_m \right),
\end{eqnarray}
the $j$-th eigenvalue of the matrix $ \mathbf{C}_i$ is called a T-eigenvalue \cite{miao2021t}, denoted by $\lambda_{i, j}$. If a symmetric T-product tensor $\mathcal{C} \in \mathbb{R}^{m \times m \times p}$ can be expressed as the format shown by Eq.~\eqref{eq:block diagonalized format}, the T-eigenvalues of  $\mathcal{C}$ with respect to the matrix $\mathbf{C}_i$ are denoted as $\lambda_{i, k_i}$, where $1 \leq k_i \leq m$, and we assume that
$\lambda_{i, 1} \geq \lambda_{i, 2} \geq \cdots \geq \lambda_{i, m}$ (including multiplicities). Then, $\lambda_{i, k_i}$ is the $k_i$-th largest T-eigenvalue associated to the matrix $\mathbf{C}_i$. If we sort all T-eigenvalues of $\mathcal{C}$ from the largest one to the smallest one, we use $\tilde{k}$, a smallest integer between 1 to $m \times p$ (inclusive) associated with $p$ given non-negative integers $k_1, k_2, \cdots, k_p$ such that there are 
$k_i$ T-eigenvalues greater or equal than $\lambda_{\tilde{k}}$ for the matrix $\mathbf{C}_{\tilde{i}}$. We set $\tilde{i}$ from $\lambda_{\tilde{k}} $ as
\begin{eqnarray}\label{eq:tilde i}
\tilde{i}= \arg \min\limits_{i }\left\{  \lambda_{\tilde{k}}= \lambda_{i, k_i} | k_i > 0 \right\}
\end{eqnarray}
Then, we will have the following Courant-Fischer theorem for T-product tensors.

\begin{theorem}\label{thm:Courant-Fischer T-product}
Given a symmetric T-product tensor $\mathcal{C} \in \mathbb{R}^{m \times m \times p}$ and $p$ non-negative  integers $k_1, k_2, \cdots, k_p$ with $0 \leq k_i \leq m$, then we have
\begin{eqnarray}
\lambda_{\tilde{k}} &=& \max\limits_{\substack{S \in \mathbb{R}^{m \times 1 \times p}\\ \dim(\mathrm{S})  = \{k_1, \cdots, k_p \}   }} \min\limits_{\mathcal{X} \in S } \frac{ \langle \mathcal{X}, \mathcal{C} \star \mathcal{X} \rangle }{ \langle \mathcal{X}, \mathcal{X} \rangle }\nonumber \\
 &=&  \min\limits_{\substack{T \in \mathbb{R}^{m \times 1 \times p}\\ \dim(T)  = \{m- k_1, \cdots, m-k_{\tilde{i}-1},  m-k_{\tilde{i}}+1,  m-k_{\tilde{i}+1},\cdots,  m-k_p \}   }} \max\limits_{\mathcal{X} \in T} \frac{ \langle \mathcal{X}, \mathcal{C} \star \mathcal{X} \rangle }{ \langle \mathcal{X}, \mathcal{X} \rangle }
\end{eqnarray}
where $\tilde{i}$ is defined by Eq.~\eqref{eq:tilde i}.
\end{theorem}
\textbf{Proof:}

First, we have to express $ \langle \mathcal{X}, \mathcal{C} \star \mathcal{X} \rangle $ by matrices of $\mathbf{C}_i$ and $\mathbf{X}_i$ through the representation shown by Eq.~\eqref{eq:block diagonalized format}. It is
\begin{eqnarray}\label{eq:thm:Courant-Fischer T-product}
\langle \mathcal{X}, \mathcal{C} \star \mathcal{X} \rangle &=& \frac{1}{p} \langle \mbox{bcirc}(\mathcal{X}), \mbox{bcirc}(\mathcal{C}) \mbox{bcirc}(\mathcal{X})    \rangle \nonumber \\
&=& \frac{1}{p} \mathrm{Tr} \left( \mbox{bcirc}(\mathcal{X})^{\mathrm{H}} \mbox{bcirc}(\mathcal{C}) \mbox{bcirc}(\mathcal{X}) \right)\nonumber \\
&=& \frac{1}{p} \mathrm{Tr} \left(\mathbf{F}^{\mathrm{H}}_p  \mbox{Diag}\left( \mathbf{x}^{\mathrm{H}}_i \mathbf{A}_i  \mathbf{x}_i: i \in \{1,\cdots,p\} \right)\mathbf{F}_p  \right)\nonumber \\
&=& \frac{1}{p} \mathrm{Tr} \left(  \mbox{Diag}\left( \mathbf{x}^{\mathrm{H}}_i \mathbf{A}_i  \mathbf{x}_i: i \in \{1,\cdots,p\} \right)  \right) = \frac{1}{p}\sum\limits_{i=1}^p \mathbf{x}^{\mathrm{H}}_i \mathbf{A}_i  \mathbf{x}_i
\end{eqnarray}

Without loss of generality, we can assume that all $k_i$ is positive since if any of these $k_i$ is zero, the term of $\langle \mathcal{X}, \mathcal{C} \star \mathcal{X} \rangle$ is reduced as $\frac{1}{p}\sum\limits_{i'} \mathbf{x}^{\mathrm{H}}_{i'} \mathbf{A}_{i'}  \mathbf{x}_{i'}$, where $k_{i'} > 0$. We will just verify the first characterization of $\lambda_{\tilde{k}}$. The other is similar. Let $S_i$ be the projection of $S$ to the space with dimension $k_i$ spanned by $\mathbf{v}_{i, 1}, \cdots, \mathbf{v}_{i, k_i}$, for every $\mathbf{x}_i \in S_{i}$, we can write $\mathbf{x}_i = \sum\limits^{k_i}_{j=1} c_{i, j}  \mathbf{v}_{i, j}$. To show that the value $\lambda_{\tilde{k}}$ is achievable, note that
\begin{eqnarray}
 \frac{ \langle \mathcal{X}, \mathcal{C} \star \mathcal{X} \rangle }{ \langle \mathcal{X}, \mathcal{X} \rangle }
&=&  \frac{  \frac{1}{p}\sum\limits_{i=1}^p \mathbf{x}^{\mathrm{H}}_i \mathbf{A}_i  \mathbf{x}_i  }{  \frac{1}{p}\sum\limits_{i=1}^p \mathbf{x}^{\mathrm{H}}_i  \mathbf{x}_i  } = \frac{  \sum\limits_{i=1}^p    \sum\limits^{k_i}_{j=1} \lambda_{i, j}    c_{i, j}^{\ast} c_{i, j}          }{   \sum\limits_{i=1}^p    \sum\limits^{k_i}_{j=1}   c_{i, j}^{\ast} c_{i, j}  } \nonumber \\
&\geq & \frac{  \sum\limits_{i=1}^p    \sum\limits^{k_i}_{j=1} \lambda_{\tilde{k}}  c_{i, j}^{\ast} c_{i, j}          }{   \sum\limits_{i=1}^p    \sum\limits^{k_i}_{j=1}   c_{i, j}^{\ast} c_{i, j}  } = \lambda_{\tilde{k}}
\end{eqnarray}
To verify that this is the maximum, let $T_{\tilde{i}}$ be the projection of $T$ to the space with dimension $k_{\tilde{i}}$ with dimension $n- k_{\tilde{i}} + 1$, then the intersection of $S$ and $T_{\tilde{i}}$ is not empty. We have
\begin{eqnarray}
\min\limits_{\mathcal{X} \in S } \frac{ \langle \mathcal{X}, \mathcal{C} \star \mathcal{X} \rangle }{ \langle \mathcal{X}, \mathcal{X} \rangle } &\leq & \min\limits_{\mathcal{X} \in S \cap T} \frac{ \langle \mathcal{X}, \mathcal{C} \star \mathcal{X} \rangle }{ \langle \mathcal{X}, \mathcal{X} \rangle }.
\end{eqnarray}
Any such $\mathbf{x}_{\tilde{i}} \in S \cap T_{\tilde{i}}$ can be expressed as $\mathbf{x}_{\tilde{i}} =  \sum\limits^{m}_{ j=k_{\tilde{i} }} c_{\tilde{i}, j}  \mathbf{v}_{\tilde{i} j}$, and any $i$ for $i \neq \tilde{i}$, we have $\mathbf{x}_{i} \in S \cap T_{i}$ expressed as $\mathbf{x}_{ i } =  \sum\limits^{m}_{ j=k_{i} + 1} c_{i, j}  \mathbf{v}_{i, j}$. Then, we have
\begin{eqnarray}
 \frac{ \langle \mathcal{X}, \mathcal{C} \star \mathcal{X} \rangle }{ \langle \mathcal{X}, \mathcal{X} \rangle }
&=&  \frac{  \frac{1}{p}\sum\limits_{i=1}^p \mathbf{x}^{\mathrm{H}}_i \mathbf{A}_i  \mathbf{x}_i  }{  \frac{1}{p}\sum\limits_{i=1}^p \mathbf{x}^{\mathrm{H}}_i  \mathbf{x}_i  } = \frac{  \sum\limits_{i=1}^p     \sum\limits^{m}_{ \substack{j=k_i + 1; i \neq \tilde{i} \\ j=k_{\tilde{i}};  i = \tilde{i}    }}   \lambda_{i, j}    c_{i, j}^{\ast} c_{i, j}          }{    \sum\limits_{i=1}^p    \sum\limits^{m}_{ \substack{j=k_i + 1; i \neq \tilde{i} \\ j=k_{\tilde{i}};  i = \tilde{i}    }   } c_{i, j}^{\ast} c_{i, j}   } \nonumber \\
&\leq & \frac{  \sum\limits_{i=1}^p    \sum\limits^{m}_{ \substack{j=k_i + 1; i \neq \tilde{i} \\ j=k_{\tilde{i}};  i = \tilde{i}    }   } \lambda_{\tilde{k}}  c_{i, j}^{\ast} c_{i, j}          }{   \sum\limits_{i=1}^p    \sum\limits^{m}_{ \substack{j=k_i + 1; i \neq \tilde{i} \\ j=k_{\tilde{i}};  i = \tilde{i}    }   } c_{i, j}^{\ast} c_{i, j}  } = \lambda_{\tilde{k}}.
\end{eqnarray}
Therefore, for all subspaces $S$ of dimensions $\{k_1, \cdots, k_p\}$, we have $\min\limits_{\mathcal{X} \in S} \frac{ \langle \mathcal{X}, \mathcal{C} \star \mathcal{X} \rangle }{ \langle \mathcal{X}, \mathcal{X} \rangle } \leq \lambda_{\tilde{k}}$
$\hfill \Box$

\section{Weyl Inequality for Symmetric T-product Tensors}\label{sec:Weyl Inequality for Symmetric T-product Tensors}

We then can apply Theorem~\ref{thm:Courant-Fischer T-product} to prove Weyl inequality for symmetric T-product tensors. 
\begin{theorem}\label{thm:Weyl Inequality for T-tensor}
Suppose $\mathcal{A}, \mathcal{B} \in \mathbb{R}^{m \times m \times p}$ are symmetric tensors with T-eigenvalues $\lambda_1 \geq \lambda_2 \geq \cdots \geq \lambda_{mp}$ and 
$\epsilon_1 \geq \epsilon_2 \geq \cdots \geq \epsilon_{mp}$, respectively. Let $\mathcal{C} = \mathcal{A} + \mathcal{B}$ with T-eigenvalues $\mu_1 \geq \mu_2 \geq \cdots \geq \mu_{mp}$. We then have:
\begin{eqnarray}
\lambda_{\tilde{k}} + \epsilon_1 \geq_1 \mu_{\tilde{k}} \geq_2 \lambda_{\tilde{k}} + \epsilon_{mp},
\end{eqnarray}
where $1 \leq \tilde{k} \leq m p$, where $\tilde{k}$ is associated to $p$ non-negative integers $k_1, \cdots, k_p$ between $0$ and $m$ inclusive.
\end{theorem}
\textbf{Proof:}
Due to Theorem~\ref{thm:Courant-Fischer T-product}, we will prove the inequality $\geq_1$ only based on  
\begin{eqnarray}
 \lambda_{\tilde{k}} &=& \min\limits_{\substack{T \in \mathbb{R}^{m \times 1 \times p}\\ \dim(T)  = \{m- k_1, \cdots, m-k_{\tilde{i}-1},  m-k_{\tilde{i}}+1,  m-k_{\tilde{i}+1},\cdots,  m-k_p \}   }} \max\limits_{\mathcal{X} \in T} \frac{ \langle \mathcal{X}, \mathcal{C} \star \mathcal{X} \rangle }{ \langle \mathcal{X}, \mathcal{X} \rangle },
\end{eqnarray}
since the inequality $\geq_2$ can be proved similarly from $ \lambda_{\tilde{k}} = \max\limits_{\substack{S \in \mathbb{R}^{m \times 1 \times p}\\ \dim(\mathrm{S})  = \{k_1, \cdots, k_p \}   }} \min\limits_{\mathcal{X} \in S } \frac{ \langle \mathcal{X}, \mathcal{C} \star \mathcal{X} \rangle }{ \langle \mathcal{X}, \mathcal{X} \rangle }$.

Because we have
\begin{eqnarray}
\lambda_{\tilde{k}} &=& \min\limits_{\substack{T \in \mathbb{R}^{m \times 1 \times p}\\ \dim(T)  = \{m- k_1, \cdots, m-k_{\tilde{i}-1},  m-k_{\tilde{i}}+1,  m-k_{\tilde{i}+1},\cdots,  m-k_p \}   }} \max\limits_{\mathcal{X} \in T} \frac{ \langle \mathcal{X}, \mathcal{C} \star \mathcal{X} \rangle }{ \langle \mathcal{X}, \mathcal{X} \rangle } \nonumber \\
& = & \min\limits_{\substack{T \in \mathbb{R}^{m \times 1 \times p}\\ \dim(T)  = \{m- k_1, \cdots, m-k_{\tilde{i}-1},  m-k_{\tilde{i}}+1,  m-k_{\tilde{i}+1},\cdots,  m-k_p \}   }} \max\limits_{\mathcal{X} \in T} \frac{ \langle \mathcal{X}, \mathcal{A} \star \mathcal{X} \rangle  +  \langle \mathcal{X}, \mathcal{B} \star \mathcal{X} \rangle             }{ \langle \mathcal{X}, \mathcal{X} \rangle }  \nonumber \\
& \leq & \min\limits_{\substack{T \in \mathbb{R}^{m \times 1 \times p}\\ \dim(T)  = \{m- k_1, \cdots, m-k_{\tilde{i}-1},  m-k_{\tilde{i}}+1,  m-k_{\tilde{i}+1},\cdots,  m-k_p \}   }} \left( \max\limits_{\mathcal{X} \in T} \frac{ \langle \mathcal{X}, \mathcal{A} \star \mathcal{X} \rangle  }{ \langle \mathcal{X}, \mathcal{X} \rangle }
+  \max\limits_{\mathcal{X} \in T} \frac{ \langle \mathcal{X}, \mathcal{B} \star \mathcal{X} \rangle  }{ \langle \mathcal{X}, \mathcal{X} \rangle } \right)  \nonumber \\
& \leq & \min\limits_{\substack{T \in \mathbb{R}^{m \times 1 \times p}\\ \dim(T)  = \{m- k_1, \cdots, m-k_{\tilde{i}-1},  m-k_{\tilde{i}}+1,  m-k_{\tilde{i}+1},\cdots,  m-k_p \}   }} \max\limits_{\mathcal{X} \in T} \frac{ \langle \mathcal{X}, \mathcal{A} \star \mathcal{X} \rangle  }{ \langle \mathcal{X}, \mathcal{X} \rangle } \nonumber \\
&  & +   \min\limits_{\substack{T \in \mathbb{R}^{m \times 1 \times p}\\ \dim(T)  = \{  \overbrace{m, \cdots, m}^{ p-\mbox{terms} } \}   } } \max\limits_{\mathcal{X} \in T} \frac{ \langle \mathcal{X}, \mathcal{B} \star \mathcal{X} \rangle  }{ \langle \mathcal{X}, \mathcal{X} \rangle }  \nonumber \\
&=&\lambda_{\tilde{k}}  + \epsilon_1.
\end{eqnarray}
Then, this theorem is proved. 
$\hfill \Box$

\section{Hanson-Wright Inequality for Random Symmetric T-product Tensors}\label{sec:Hanson-Wright Inequality for Random Symmetric T-product Tensors}

For random variables, Bernstein inequalities give the upper tail of a sum of
independent, zero-mean random variables that are either bounded or subexponential. In Theorem 1.7 at~\cite{chang2021t_P_II}, we proved Bernstein bounds for a sum of zero-mean random T-product tensors. 

\begin{theorem}[T-product Tensor Bernstein Bounds with Bounded $\lambda_{\max}$]\label{thm:Bounded Tensor Bernstein}
Given a finite sequence of independent Hermitian T-product tensors $\{ \mathcal{X}_i  \in \mathbb{C}^{m \times m \times p} \}$ that satisfy
\begin{eqnarray}\label{eq1:thm:Bounded Tensor Bernstein}
\mathbb{E} \mathcal{X}_i = 0 \mbox{~~and~~} \lambda_{\max}(\mathcal{X}_i) \leq T 
\mbox{~~almost surely.} 
\end{eqnarray}
Define the total varaince $\sigma^2$ as: $\sigma^2 \define \left\Vert \sum\limits_{i=1}^n \mathbb{E} \left( \mathcal{X}^2_i \right) \right\Vert$.
Then, we have following inequalities:
\begin{eqnarray}\label{eq2:thm:Bounded Tensor Bernstein}
\mathrm{Pr} \left( \lambda_{\max}\left( \sum\limits_{i=1}^{n} \mathcal{X}_i \right)\geq \theta \right) \leq mp \exp \left( \frac{-\theta^2/2}{\sigma^2 + T\theta/3}\right);
\end{eqnarray}
and
\begin{eqnarray}\label{eq3:thm:Bounded Tensor Bernstein}
\mathrm{Pr} \left( \lambda_{\max}\left( \sum\limits_{i=1}^{n} \mathcal{X}_i \right)\geq \theta \right) \leq mp \exp \left( \frac{-3 \theta^2}{ 8 \sigma^2}\right)~~\mbox{for $\theta \leq \sigma^2/T$};
\end{eqnarray}
and
\begin{eqnarray}\label{eq4:thm:Bounded Tensor Bernstein}
\mathrm{Pr} \left( \lambda_{\max}\left( \sum\limits_{i=1}^{n} \mathcal{X}_i \right)\geq \theta \right) \leq mp \exp \left( \frac{-3 \theta}{ 8 T } \right)~~\mbox{for $\theta \geq \sigma^2/T$}.
\end{eqnarray}
\end{theorem}

Finally, we present the Hanson-Wright Inequality for random symmetric T-product tensors.
\begin{theorem}[Hanson-Wright Inequality for Random Symmetric T-product Tensors]{thm}{HWThm}\label{thm:T-HW inequality}
We define a vector of random T-product tensors $\overline{\mathcal{X}} \in \mathbb{R}^{(n \times m)  \times m \times p}$ as:
\begin{eqnarray}\label{eq:vec X def 1}
\overline{\mathcal{X}} = \begin{bmatrix}
           \mathcal{X}_{1} \\
           \mathcal{X}_{2} \\
           \vdots \\
           \mathcal{X}_{n}
         \end{bmatrix},
\end{eqnarray}
where random symmetric T-product tensors $\mathcal{X}_{i} \in \mathbb{R}^{m \times m \times p} $ are independent random symmetric T-product tensors with $\mathbb{E} \mathcal{X}_i = \mathcal{O}$ for $1 \leq i \leq n$. We also require another fixed tensor $\overline{\overline{\mathcal{A}}} \in  \mathbb{R}^{(n \times m) \times (n \times m) \times p}$, which is defined as:
\begin{eqnarray}\label{eq:matrix A def 1}
\overline{\overline{\mathcal{A}}} = \begin{bmatrix}
           \mathcal{A}_{1,1} &  \mathcal{A}_{1,2} & \cdots & \mathcal{A}_{1, n}  \\
           \mathcal{A}_{2,1} &  \mathcal{A}_{2,2} & \cdots & \mathcal{A}_{2, n}  \\
           \vdots & \vdots &   \vdots & \vdots \\
           \mathcal{A}_{n,1} &  \mathcal{A}_{n,2} & \cdots & \mathcal{A}_{n, n}  \\
         \end{bmatrix},
\end{eqnarray}
where $\mathcal{A}_{i,j} \in \mathbb{R}^{m \times m \times p}$ are symmetric T-product tensors also. We also require following assumptions. Define random Hermitian tensor $\mathcal{Y}_i$ as 
\begin{eqnarray}\label{eq:diag Yi def main}
\mathcal{Y}_i \define \mathcal{X}_i \star \mathcal{A}_{i,i} \star \mathcal{X}_i  - \mathbb{E}\left( \mathcal{X}_i \star \mathcal{A}_{i,i} \star  \mathcal{X}_i \right),~~\mbox{for $1 \leq i \leq n$;}
\end{eqnarray}
we assume that 
\begin{eqnarray}
\mathbb{E} \mathcal{Y}_i = \mathcal{O} \mbox{~~and~~} \lambda_{\max}(\mathcal{Y}_i) \leq T_{\mbox{dg}} 
\mbox{~~almost surely.} 
\end{eqnarray}
Define the total varaince $\sigma_{\mbox{dg}}^2$ as: $\sigma_{\mbox{dg}}^2 \define \left\Vert \sum\limits_{i=1}^n \mathbb{E} \left( \mathcal{Y}^2_i \right) \right\Vert$, where $\left\Vert \cdot \right\Vert$ represents the spectral norm, which equals the largest singular value of a T-product tensor.

Moreover, we define random Hermitian tensor $\mathcal{Z}_k$ for $k=1,2,\cdots, n^2 - n$ as 
\begin{eqnarray}\label{eq:cp Zk def main}
\mathcal{Z}_k \define  \mathcal{X}^{(1)}_i \star  \mathcal{A}_{i,j} \star  \mathcal{X}^{(2)}_j~~\mbox{for $1 \leq i \neq j \leq n$;}
\end{eqnarray}
where the tensors $\mathcal{X}^{(1)}_i$ are identical distribution copy for the tensors $\mathcal{X}_i$, and the tensors $\mathcal{X}^{(2)}_j$ are identical distribution copy for the tensors $\mathcal{X}_j$, then we assume that 
\begin{eqnarray}
\mathbb{E} \mathcal{Z}_k=\mathcal{O} 
\mbox{~~almost surely.} 
\end{eqnarray}
Given any realization of the random tensor $\mathcal{X}_i$, denoted as $\tilde{\mathcal{X}}_i$, we assume that 
\begin{eqnarray}\label{eq:max e-value bound cp TP}
\lambda_{\max}\left( \tilde{\mathcal{X}}_i \star_M \mathcal{A}_{i, j} \star_M \mathcal{X}_{j} \right) \leq T_{\mbox{cp}} \mbox{~~almost surely},
\end{eqnarray}
where $T_{\mbox{cp}}$ is a positive real number, and all $1 \leq i \neq j \leq n$. The total variance with respect to the tensor $\tilde{\mathcal{X}}_i$ is defined as
\begin{eqnarray}\label{eq:sigma cp TP}
\sigma^2_{\mbox{cp}}( \tilde{\mathcal{X}}_i)  \define  \left\Vert \sum\limits_{j=1, \neq i}^{n} \mathbb{E} \left( \left(\tilde{\mathcal{X}}_i \star_M \mathcal{A}_{i, j} \star_M \mathcal{X}_{j}   \right)^2 \right) \right\Vert.
\end{eqnarray}
The function $f(\tilde{\mathcal{X}}_i)$ is the probability density function for the realization tensor $\tilde{\mathcal{X}}_i$.

Then, we have
\begin{eqnarray}
\mathrm{Pr}\left(  \lambda_{\max} \left( \overline{\mathcal{X}}^{\mathrm{T}} 
\overline{\overline{\mathcal{A}}} \overline{\mathcal{X}} - \mathbb{E}\left( \overline{\mathcal{X}}^{\mathrm{T}} 
\overline{\overline{\mathcal{A}}} \overline{\mathcal{X}} \right)  \right) \geq \theta   \right) \leq  
\overbrace{\mathrm{Pr}\left( \lambda_{\max}\left(  \sum\limits_{1 \leq i \neq j \leq n}\mathcal{X}_i \star \mathcal{A}_{i,j} \star  \mathcal{X}_j \right) \geq \frac{\theta}{2}  \right) }^{\define \mathrm{P}_{\mbox{cp}} } \nonumber \\
 + \overbrace{\mathrm{Pr} \left( \lambda_{\max}\left(  \sum\limits_{i=1}^{n} \left(  \mathcal{X}_i  \star \mathcal{A}_{i,i} \star \mathcal{X}_i - \mathbb{E} \left(  \mathcal{X}_i  \star \mathcal{A}_{i,i} \star \mathcal{X}_i \right)    \right) \right) \geq \frac{\theta}{2} \right) }^{\define  \mathrm{P}_{\mbox{dg}}} ~~~~~~~~~~~~~~~~~~~~~~~~~~~~~~~~~~ \nonumber \\
 \leq   C_4 mp \sum\limits_{i=1}^{n}  \int_{\tilde{\mathcal{X}}_i}  \exp \left( \frac{ - \theta^2}{  8 n^2C^2_4 \sigma^2_{\mbox{cp}}( \tilde{\mathcal{X}}_i) + 4T_{\mbox{cp}}\theta n C_4 /3 }\right)  f(\tilde{\mathcal{X}}_i) d \tilde{\mathcal{X}}_i \nonumber \\
+  mp \exp \left( \frac{-\theta^2}{8 \sigma_{\mbox{dg}}^2 + 4 T_{\mbox{dg}}\theta/3}\right),~~~~~~~~~~~~~~~~~~~~~~~~~~~~~~~~~~
\end{eqnarray}
where $\mathrm{P}_{\mbox{cp}}$ and $\mathrm{P}_{\mbox{dg}}$ are probability bounds related to the coupling sum and the diagonal sum parts, respectively, and the term $C_4$ is a positive constant. 

If  $\frac{\theta}{2n C_4} \leq \frac{\sigma^2_{\mbox{cp}}( \tilde{\mathcal{X}}_i)  }{T_{\mbox{cp}}}$ with respect to $i=1,2,\ldots,n$ and $\theta \leq 2 \sigma_{\mbox{dg}}^2 / T_{\mbox{dg}}$, we have 
\begin{eqnarray}
\mathrm{Pr}\left(  \lambda_{\max} \left( \overline{\mathcal{X}}^{\mathrm{T}} 
\overline{\overline{\mathcal{A}}} \overline{\mathcal{X}} - \mathbb{E}\left( \overline{\mathcal{X}}^{\mathrm{T}} 
\overline{\overline{\mathcal{A}}} \overline{\mathcal{X}} \right)  \right) \geq \theta   \right) &\leq&  
 mp C_4 \sum\limits_{i=1}^{n}   \int_{\tilde{\mathcal{X}}_i}  \exp \left( \frac{ - 3 \theta^2}{   32 n^2C^2_4 \sigma^2_{\mbox{cp}}( \tilde{\mathcal{X}}_i) }\right)  f(\tilde{\mathcal{X}}_i) d \tilde{\mathcal{X}}_i  \nonumber \\
&  & + mp \exp \left( \frac{-3 \theta^2}{ 32 \sigma_{\mbox{dg}}^2}\right). 
\end{eqnarray}
Moreover, if $\frac{\theta}{2n C_4} \geq \frac{\sigma^2_{\mbox{cp}}( \tilde{\mathcal{X}}_i)  }{T_{\mbox{cp}}}$ with respect to $i=1,2,\ldots,n$, and $\theta \geq 2 \sigma_{\mbox{dg}}^2 / T_{\mbox{dg}}$, we have 
\begin{eqnarray}
\mathrm{Pr}\left(  \lambda_{\max} \left( \overline{\mathcal{X}}^{\mathrm{T}} 
\overline{\overline{\mathcal{A}}} \overline{\mathcal{X}} - \mathbb{E}\left( \overline{\mathcal{X}}^{\mathrm{T}} 
\overline{\overline{\mathcal{A}}} \overline{\mathcal{X}} \right)  \right) \geq \theta   \right) &\leq&  
mnp C_4  \exp \left( \frac{ - 3 \theta}{  16 nC_4 T_{\mbox{cp}}         }\right)  + mp \exp \left( \frac{-3 \theta}{ 16 T_{\mbox{dg}} } \right). 
\end{eqnarray}
\end{theorem}
\textbf{Proof:}
Since the proof arguments of Theorem~\ref{thm:HW inequality} can still be valid for T-product tensors by applying Theorem~\ref{thm:Weyl Inequality for T-tensor}, this theorem is proved by using Theorem~\ref{thm:Bounded Tensor Bernstein} to modify Lemma~\ref{lma:diag sum boud} and Lemma~\ref{lma:cp sum boud} for tensors under T-product. 
$\hfill \Box$

\end{appendices}

\bibliographystyle{IEEETran}
\bibliography{EProd_QuadSum_Bib}

\begin{thebibliography}{10}
\providecommand{\url}[1]{#1}
\csname url@samestyle\endcsname
\providecommand{\newblock}{\relax}
\providecommand{\bibinfo}[2]{#2}
\providecommand{\BIBentrySTDinterwordspacing}{\spaceskip=0pt\relax}
\providecommand{\BIBentryALTinterwordstretchfactor}{4}
\providecommand{\BIBentryALTinterwordspacing}{\spaceskip=\fontdimen2\font plus
\BIBentryALTinterwordstretchfactor\fontdimen3\font minus
  \fontdimen4\font\relax}
\providecommand{\BIBforeignlanguage}[2]{{%
\expandafter\ifx\csname l@#1\endcsname\relax
\typeout{** WARNING: IEEEtran.bst: No hyphenation pattern has been}%
\typeout{** loaded for the language `#1'. Using the pattern for}%
\typeout{** the default language instead.}%
\else
\language=\csname l@#1\endcsname
\fi
#2}}
\providecommand{\BIBdecl}{\relax}
\BIBdecl

\bibitem{adamczak2015note}
R.~Adamczak, ``A note on the hanson-wright inequality for random vectors with
  dependencies,'' \emph{Electronic Communications in Probability}, vol.~20, pp.
  1--13, 2015.

\bibitem{hanson1971bound}
D.~L. Hanson and F.~T. Wright, ``A bound on tail probabilities for quadratic
  forms in independent random variables,'' \emph{The Annals of Mathematical
  Statistics}, vol.~42, no.~3, pp. 1079--1083, 1971.

\bibitem{vershynin2018high}
R.~Vershynin, \emph{High-dimensional probability: An introduction with
  applications in data science}.\hskip 1em plus 0.5em minus 0.4em\relax
  Cambridge university press, 2018, vol.~47.

\bibitem{krahmer2014suprema}
F.~Krahmer, S.~Mendelson, and H.~Rauhut, ``Suprema of chaos processes and the
  restricted isometry property,'' \emph{Communications on Pure and Applied
  Mathematics}, vol.~67, no.~11, pp. 1877--1904, 2014.

\bibitem{qi2017tensor}
L.~Qi and Z.~Luo, \emph{Tensor analysis: spectral theory and special
  tensors}.\hskip 1em plus 0.5em minus 0.4em\relax SIAM, 2017.

\bibitem{wu2010robust}
Q.~Wu, L.~Zhang, and G.~Shi, ``Robust multifactor speech feature extraction
  based on gabor analysis,'' \emph{IEEE Transactions on Audio, Speech, and
  Language Processing}, vol.~19, no.~4, pp. 927--936, Aug. 2010.

\bibitem{mirsamadi2016generalized}
S.~Mirsamadi and J.~H. Hansen, ``A generalized nonnegative tensor factorization
  approach for distant speech recognition with distributed microphones,''
  \emph{IEEE/ACM Transactions on Audio, Speech, and Language Processing},
  vol.~24, no.~10, pp. 1721--1731, Jun. 2016.

\bibitem{muti2007survey}
D.~Muti and S.~Bourennane, ``Survey on tensor signal algebraic filtering,''
  \emph{Signal Processing}, vol.~87, no.~2, pp. 237--249, Feb. 2007.

\bibitem{shen2020topology}
Y.~Shen, X.~Fu, G.~B. Giannakis, and N.~D. Sidiropoulos, ``Topology
  identification of directed graphs via joint diagonalization of correlation
  matrices,'' \emph{IEEE Transactions on Signal and Information Processing over
  Networks}, vol.~6, pp. 271--283, Apr. 2020.

\bibitem{shen2017tensor}
Y.~Shen, B.~Baingana, and G.~B. Giannakis, ``{T}ensor {D}ecompositions for
  {I}dentifying {D}irected graph {T}opologies and {T}racking {D}namic
  {N}etworks,'' \emph{IEEE Transactions on Signal Processing}, vol.~65, no.~14,
  pp. 3675--3687, Apr. 2017.

\bibitem{fu2015joint}
X.~Fu, K.~Huang, W.-K. Ma, N.~D. Sidiropoulos, and R.~Bro, ``Joint tensor
  factorization and outlying slab suppression with applications,'' \emph{IEEE
  Transactions on Signal Processing}, vol.~63, no.~23, pp. 6315--6328, Aug.
  2015.

\bibitem{ko2020fast}
C.-Y. Ko, K.~Batselier, L.~Daniel, W.~Yu, and N.~Wong, ``Fast and accurate
  tensor completion with total variation regularized tensor trains,''
  \emph{IEEE Transactions on Image Processing}, May 2020.

\bibitem{jiang2020framelet}
T.-X. Jiang, M.~K. Ng, X.-L. Zhao, and T.-Z. Huang, ``Framelet representation
  of tensor nuclear norm for third-order tensor completion,'' \emph{IEEE
  Transactions on Image Processing}, vol.~29, pp. 7233--7244, Jun. 2020.

\bibitem{de2008constrained}
A.~L. de~Almeida, G.~Favier, and J.~C.~M. Mota, ``Constrained tensor modeling
  approach to blind multiple-antenna {CDMA} schemes,'' \emph{IEEE Transactions
  on Signal Processing}, vol.~56, no.~6, pp. 2417--2428, May 2008.

\bibitem{zhijin2018blind}
Z.~Zhijin, Y.~Hui, and J.~S. Fangfang~QIANG, ``Blind estimation of spreading
  codes for multi-antenna lc-ds-cdma signals based on tensor decomposition,''
  \emph{Journal on Communications}, vol.~39, no.~10, p.~52, 2018.

\bibitem{nion2010tensor}
D.~Nion and N.~D. Sidiropoulos, ``Tensor algebra and multidimensional harmonic
  retrieval in signal processing for {MIMO} radar,'' \emph{IEEE Transactions on
  Signal Processing}, vol.~58, no.~11, pp. 5693--5705, Jul. 2010.

\bibitem{sidiropoulos2000parallel}
N.~D. Sidiropoulos, R.~Bro, and G.~B. Giannakis, ``Parallel factor analysis in
  sensor array rrocessing,'' \emph{IEEE Transactions on Signal Processing},
  vol.~48, no.~8, pp. 2377--2388, Aug. 2000.

\bibitem{wang2019neural}
X.~Wang, M.~Che, and Y.~Wei, ``Neural networks based approach solving
  multi-linear systems with m-tensors,'' \emph{Neurocomputing}, vol. 351, pp.
  33--42, 2019.

\bibitem{MR3395816}
\BIBentryALTinterwordspacing
W.~Ding, L.~Qi, and Y.~Wei, ``Fast {H}ankel tensor-vector product and its
  application to exponential data fitting,'' \emph{Numer. Linear Algebra
  Appl.}, vol.~22, no.~5, pp. 814--832, 2015. [Online]. Available:
  \url{https://doi.org/10.1002/nla.1970}
\BIBentrySTDinterwordspacing

\bibitem{MR3947912}
\BIBentryALTinterwordspacing
H.-R. Xu, D.-H. Li, and S.-L. Xie, ``An equivalent tensor equation to the
  tensor complementarity problem with positive semi-definite {$Z$}-tensor,''
  \emph{Optim. Lett.}, vol.~13, no.~4, pp. 685--694, 2019. [Online]. Available:
  \url{https://doi.org/10.1007/s11590-018-1268-4}
\BIBentrySTDinterwordspacing

\bibitem{MR3479021}
\BIBentryALTinterwordspacing
L.-B. Cui, C.~Chen, W.~Li, and M.~K. Ng, ``An eigenvalue problem for even order
  tensors with its applications,'' \emph{Linear Multilinear Algebra}, vol.~64,
  no.~4, pp. 602--621, 2016. [Online]. Available:
  \url{https://doi.org/10.1080/03081087.2015.1071311}
\BIBentrySTDinterwordspacing

\bibitem{anandkumar2015tensor}
A.~Anandkumar, R.~Ge, D.~Hsu, S.~M. Kakade, and M.~Telgarsky, ``Tensor
  decompositions for learning latent variable models ({A} survey for {ALT}),''
  in \emph{Proceedings of International Conference on Algorithmic Learning
  Theory}.\hskip 1em plus 0.5em minus 0.4em\relax Springer, Oct. 2015, pp.
  19--38.

\bibitem{sidiropoulos2017tensor}
N.~D. Sidiropoulos, L.~De~Lathauwer, X.~Fu, K.~Huang, E.~E. Papalexakis, and
  C.~Faloutsos, ``Tensor decomposition for signal processing and machine
  learning,'' \emph{IEEE Transactions on Signal Processing}, vol.~65, no.~13,
  pp. 3551--3582, Jul. 2017.

\bibitem{MR3616422}
R.~Gurau, \emph{Random tensors}.\hskip 1em plus 0.5em minus 0.4em\relax Oxford
  University Press, Oxford, 2017.

\bibitem{MR3783911}
\BIBentryALTinterwordspacing
I.~R. Klebanov and G.~Tarnopolsky, ``Uncolored random tensors, melon diagrams,
  and the {S}achdev-{Y}e-{K}itaev models,'' \emph{Phys. Rev. D}, vol.~95,
  no.~4, pp. 046\,004, 13, 2017. [Online]. Available:
  \url{https://doi.org/10.1103/physrevd.95.046004}
\BIBentrySTDinterwordspacing

\bibitem{MR4140540}
\BIBentryALTinterwordspacing
R.~Vershynin, ``Concentration inequalities for random tensors,''
  \emph{Bernoulli}, vol.~26, no.~4, pp. 3139--3162, 2020. [Online]. Available:
  \url{https://doi.org/10.3150/20-BEJ1218}
\BIBentrySTDinterwordspacing

\bibitem{chang2021tensor}
S.~Y. Chang, ``Tensor expander chernoff bounds,'' 2021.

\bibitem{chang2021TProdI}
------, ``T product tensors part i: Inequalities,'' \emph{arXiv preprint
  arXiv:2107.06285}, 2021.

\bibitem{chang2021TProdII}
------, ``T product tensors part ii: Tail bounds for sums of random t product
  tensors,'' \emph{arXiv preprint arXiv:2107.06224}, 2021.

\bibitem{chang2020convenient}
------, ``Convenient tail bounds for sums of random tensors,'' 2020.

\bibitem{chang2021general}
------, ``General tail bounds for random tensors summation: Majorization
  approach,'' 2021.

\bibitem{HW_T_SYChang_2021}
------, ``{H}anson-{W}right inequality for random tensors under {T}-product,''
  2021.

\bibitem{MR3913666}
\BIBentryALTinterwordspacing
M.~Liang and B.~Zheng, ``Further results on {M}oore-{P}enrose inverses of
  tensors with application to tensor nearness problems,'' \emph{Comput. Math.
  Appl.}, vol.~77, no.~5, pp. 1282--1293, 2019. [Online]. Available:
  \url{https://doi.org/10.1016/j.camwa.2018.11.001}
\BIBentrySTDinterwordspacing

\bibitem{ni2019hermitian}
G.~Ni, ``Hermitian tensor and quantum mixed state,'' \emph{arXiv preprint
  arXiv:1902.02640}, 2019.

\bibitem{de1993bounds}
V.~H. de~la Pe{\~n}a and S.~J. Montgomery-Smith, ``Bounds on the tail
  probability of {U}-statistics and quadratic forms,'' \emph{arXiv preprint
  math/9309210}, 1993.

\bibitem{kilmer2013third}
M.~E. Kilmer, K.~Braman, N.~Hao, and R.~C. Hoover, ``Third-order tensors as
  operators on matrices: A theoretical and computational framework with
  applications in imaging,'' \emph{SIAM Journal on Matrix Analysis and
  Applications}, vol.~34, no.~1, pp. 148--172, 2013.

\bibitem{miao2021t}
Y.~Miao, L.~Qi, and Y.~Wei, ``{T}-{J}ordan canonical form and {T}-{D}razin
  inverse based on the {T}-product,'' \emph{Communications on Applied
  Mathematics and Computation}, vol.~3, no.~2, pp. 201--220, 2021.

\bibitem{chang2021t_P_II}
S.~Y. Chang, ``T product tensors part {II}: Tail bounds for sums of random
  {T}-product tensors,'' 2021.

\end{thebibliography}

\end{document}